\documentclass[12pt,reqno]{amsart}
\usepackage{lmodern}
\usepackage[T1]{fontenc}

\usepackage[numbers]{natbib}
\usepackage{amsmath,amssymb,latexsym} 
\usepackage{fullpage}
\usepackage{etoolbox}
\usepackage{enumitem}
\usepackage{color}
\usepackage[colorlinks,linkcolor=blue,citecolor=magenta]{hyperref}
\usepackage[colorinlistoftodos,bordercolor=orange,backgroundcolor=orange!20,linecolor=orange,textsize=scriptsize]{todonotes}
\usepackage{soul}
\usepackage{microtype}
\usepackage{yhmath}
\usepackage{float}
\usepackage{algorithm}
\usepackage[noend]{algpseudocode}
\usepackage{amsfonts}
\usepackage{booktabs}
\usepackage{comment}
\usepackage{mathtools}
\usepackage{calc}
\usepackage{float}
\usepackage{pstricks}
\usepackage{pstricks-add}
\usepackage{eurosym}
\usepackage{cleveref}
\usepackage{pgfplots}
\usepackage{multirow}

\usepackage{geometry}
\definecolor{CompactWColor}{RGB}{0,114,178}   
\definecolor{CompactSColor}{RGB}{102,194,255} 
\definecolor{CGWColor}{RGB}{213,94,0}         
\definecolor{CGSColor}{RGB}{230,159,0}        

\theoremstyle{definition}

\DeclareFontFamily{U}{mathx}{}
\DeclareFontShape{U}{mathx}{m}{n}{ <-> mathx10 }{}
\DeclareSymbolFont{mathx}{U}{mathx}{m}{n}
\DeclareFontSubstitution{U}{mathx}{m}{n}
\DeclareMathAccent{\widecheck}{0}{mathx}{"71}

\title{A column-generation approach for an electricity technician routing and scheduling problem with a lexicographic objective}

\author{Elise Bangerter}
\address{E. Bangerter, Decision Support \& Operations Research Group, Department of Informatics, University of Fribourg, Fribourg, Switzerland}
\email{elise.bangerter@unifr.ch}
\author{David Schindl}
\address{D. Schindl, University of Applied Sciences Western Switzerland, Geneva, Switzerland\\
Decision Support \& Operations Research Group, Department of Informatics, University of Fribourg, Fribourg, Switzerland}
\email{david.schindl@hesge.ch}
\author{Meritxell Pacheco Paneque}
\address{M. Pacheco Paneque, Decision Support \& Operations Research Group, Department of Informatics, University of Fribourg, Fribourg, Switzerland}
\email{meritxell.pacheco@unifr.ch}
\author{Nour Elhouda Tellache}
\address{N.E.H. Tellache, Decision Support \& Operations Research Group, Department of Informatics, University of Fribourg, Fribourg, Switzerland}
\email{nourelhouda.tellache@unifr.ch}
\author{Rodolphe Griset}
\address{R. Griset, Hopia, Paris, 75009, France \\ Former Affiliation: EDF R\&D, OSIRIS Department, Palaiseau, 91120, France}
\email{rodolphe.griset@hotmail.fr}


\begin{document}
	
	\maketitle
	
	\begin{abstract}
	Electric utility companies perform numerous technical interventions every day. Since it is generally not possible to complete all planned interventions within a single day, companies face two objectives: maximizing the total duration of completed interventions (primary objective) and minimizing the associated operational cost (secondary objective). In this paper, we introduce a multi-objective variant of the technician routing and scheduling problem in which both objectives are optimized in lexicographic order. We propose a compact mixed-integer linear formulation and an extended set-packing-based formulation. To handle the objectives within a single-objective framework, we consider weighted-sum reformulations that preserve lexicographic priorities as well as sequential reformulations that individually optimize each objective while maintaining the optimal value of higher-priority ones. For the extended formulation, we develop an exact column-generation-based algorithm, in which the pricing subproblems are solved via a labeling algorithm based on dynamic programming. As technician schedules are typically generated on a daily basis, the algorithm is designed to deliver high-quality solutions within short computation times (e.g., 5 minutes). Computational experiments on real-life instances provided by the French electric utility company show that the CG-based algorithm proves optimality on a larger number of small instances than the compact formulation and consistently outperforms it on larger instances. In particular, the sequential CG-based variant finds the best-known solutions on more instances and achieves lower mean gaps relative to the best solution found in each instance category.
	\end{abstract}

\section{Introduction}
\label{sec:introduction}

The operation and maintenance of electricity distribution networks involve numerous technical interventions performed daily. Electric utility companies aim to execute these interventions efficiently while minimizing the associated operational costs. The problem addressed in this paper arises in this context and originates from a technician routing and scheduling problem (TRSP) at Électricité de France (EDF), the French electric utility company.

Specifically, we consider a crew of technicians and a set of interventions that might be performed in a given day. Technicians are pre-assigned into fixed teams that stay together in the same vehicle for the entire day. Each intervention has a time window and a set of required skills, and each team has an associated skill profile. The problem is to assign teams to interventions and dispatch them on routes that satisfy skill requirements, respect time windows, and comply with additional operational constraints such as a mandatory lunch break and workload limits. Since it is generally not possible to complete all interventions within a single day, minimizing total operational cost alone would not be a meaningful objective. Indeed, if the only goal were to minimize costs, and no constraints required any interventions to be performed, the optimal solution would be to perform none, resulting in zero cost. Unperformed interventions result in unmet contractual obligations or delayed service, so maximizing coverage is essential. Therefore, the problem features two objectives: the primary objective is to maximize the total duration of covered interventions, and the secondary objective is to minimize total operational cost, including personnel and distance-based travel expenses.

Multi-objective optimization methods typically aim to explore trade-offs between conflicting objectives, often represented by the Pareto front. These methods are well-suited when objectives are considered equally important or when a balance between them is desired. In contrast, lexicographic optimization imposes a strict priority order among objectives. Rather than seeking compromise, the primary objective is optimized first, and secondary objectives are only considered when multiple solutions yield the same value for higher-priority objectives. This structure induces a total order over the solution space, meaning that the concept of a Pareto front does not apply.

This problem can be seen as a variant of the TRSP, distinguished primarily by its lexicographic objective function. To the best of our knowledge, the TRSP with lexicographically-ordered objectives has not yet been addressed in the literature. It also incorporates operational constraints specific to EDF, particularly technician workload limitations to deal with unexpected events, and allows certain interventions, classified as long, to be interrupted for lunch breaks.

In this paper, we propose two mathematical formulations for the problem: a compact formulation and a set-packing formulation. We consider two strategies to handle the lexicographic objective function within a single-objective framework: a weighted reformulation that enforces the lexicographic hierarchy of objectives, and a sequential reformulation in which objectives are optimized one at a time subject to preserving previously attained optimal values. For the set-packing formulation, we develop a column-generation (CG)–based algorithm that is adapted to both reformulations. While our approach includes heuristic components to guide the search, it remains an exact method based on~CG.

The proposed reformulations and CG-based algorithm variants are evaluated through an extensive computational study based on real operational data provided by EDF. The experiments assess both solution quality and computational efficiency under time constraints representative of practical settings, with a particular emphasis on short runtimes required for daily planning. The performance of the CG–based algorithm is compared with that of the compact formulation across a range of instance sizes, including standardized test sets and large real-life instances, allowing us to analyze scalability and robustness under realistic operational conditions.

To sum up, the main contributions of this paper are threefold. First, we introduce a multi-objective variant of the TRSP with a strict lexicographic ordering of objectives, in which intervention coverage takes priority over operational cost. Second, we propose two mixed-integer linear programming (MILP) formulations (compact and set-packing) and develop two single-objective reformulations for each (weighted and sequential) that preserve the lexicographic hierarchy. Third, we design an exact CG-based algorithm for the set-packing formulation, which can also be applied heuristically for large instances under a restricted time budget. In addition, all proposed solution methodologies are evaluated on a set of instances derived from real-life operational data provided by EDF that have been standardized to enable meaningful computational comparisons.

The remainder of the paper is organized as follows. Section~\ref{sec:LR} reviews the related literature. Section~\ref{sec:problem_def} formally defines the problem. Section~\ref{sec:math_form} presents the compact and set-packing formulations. Sections~\ref{sec:solution_methods} and~\ref{sec:ColGen} describe the proposed single-objective reformulations and the CG-based algorithm, respectively. Section~\ref{sec:comp_result} reports the computational experiments and discusses managerial insights derived from the obtained solutions. Finally, Section~\ref{sec:conclusion} concludes the paper and outlines directions for future work.

\section{Literature review}\label{sec:LR}
The TRSP is a specific application within the broader class of workforce scheduling and routing problems (WSRP). WSRPs involve the scheduling and routing of skilled personnel to carry out tasks at different locations. Examples include home healthcare, where nurses visit patients to provide treatment; security services, where guards perform patrols across different sites; and technical field services, where technicians handle installations, maintenance, or repairs. We refer the reader to~\cite{castillo2016workforce} for a review of WSRPs.

The earliest work to our knowledge that addresses the TRSP is presented by~\cite{tsang1997fast}, in which a set of tasks must be performed in a day by a team of engineers at British Telecom. Each task is defined by its location, duration, and temporal constraints, indicating whether it must be scheduled in the morning, afternoon, as the first or last task of the day, or at any time. Engineers are characterized by their base location, working hours (start and end times), overtime limits, and a skill level represented by a coefficient between 0 and 1, indicating the proportion of the task duration required based on their expertise. The objective is to schedule and route engineers such that all constraints are satisfied while minimizing the total operational cost. To address this problem, the authors propose a local search-based heuristic. A more general variant is studied in~\cite{kovacs2012adaptive}, where technicians are characterized by multiple skills, each with an associated proficiency level. Tasks require specific skills at or above certain minimum levels, and may be outsourced if no available technician meets the requirements. The authors consider two problem settings: one involving team formation, in which technicians are grouped into fixed teams that remain together throughout the day; and another without team formation, where each technician operates independently. The objective is to minimize the total cost, including both travel and outsourcing expenses. To solve the problem, the authors use an adaptive large neighborhood search (ALNS) algorithm. Pillac et al.~\cite{pillac2013parallel} address a more constrained variant incorporating resource management into the problem. In their formulation, each technician is equipped with a set of tools and spare parts, while each task requires a specific subset of these items. Tools are considered as renewable resources, whereas spare parts are considered non-renewable. Technicians start their routes from home with an initial stock of resources and may replenish them during the day at a central depot. To tackle this setting, the authors propose a matheuristic approach that integrates a constructive heuristic, a parallel ALNS, and mathematical programming.

Mathlouthi et al.\cite{mathlouthi2018mixed, mathlouthi2021branch} consider a variant of the TRSP, referred to as the multi-attribute TRSP, which integrates resource management for spare and special parts, as well as multiple time windows for each task. The objective is to maximize total profit, defined as the revenue from completed tasks minus operational costs related to travel and technician overtime. The authors propose a MILP model in~\cite{mathlouthi2018mixed}, and a branch-and-price algorithm in~\cite{mathlouthi2021branch}, where the master problem is formulated as a set packing model and the pricing subproblem as an elementary shortest path problem with resource constraints (ESPPRC). The pricing subproblem is solved using the algorithm proposed by~\cite{feillet2004exact}, which is based on the Desrochers’ label correcting algorithm~\cite{desrochers1988algorithm}.

In~\cite{tricoire2013exact}, the authors investigate a multi-period, multi-depot variant of the TRSP, where the planning horizon extends over several days, and task time windows may also span multiple days. Technicians are allowed to start and end their routes either from their homes or from a central depot, with the possibility of differing start and end locations. Additionally, technician availability can vary across days. The objective is to minimize the total distance traveled. To solve the problem, the authors propose a branch-and-price algorithm, in which the pricing subproblem is modeled as an ESPPRC. They employ both exact and heuristic methods for solving the ESPPRC. The exact method adapts the dynamic programming algorithm of~\cite{feillet2004exact}, which in turn builds on the work of~\cite{desrochers1992new}. In~\cite{zamorano2017branch}, the authors also address a multi-period TRSP, extending the problem to include daily team formation in addition to technician-task assignment and routing decisions. The objective is to minimize total operational costs, which include travel expenses, customer waiting times, and technician overtime. The authors propose a mixed-integer programming model and a branch-and-price algorithm, exploring two different decomposition strategies. The pricing subproblem is also formulated as an ESPPRC and solved using a dynamic programming approach similar to that of~\cite{feillet2004exact}.

The authors in~\cite{cortes2014branch} address a variant of the TRSP that incorporates soft time windows. The objective is to minimize a total cost function that accounts for time-window violations, total travel time, and penalties for unserved tasks. To solve the problem, they propose a branch-and-price algorithm, in which the pricing subproblem is handled using constraint programming. In~\cite{pourjavad2022optimization}, the authors investigate an overnight variant of the TRSP, where technicians are allowed to perform multi-day tours, spending the night in one of the serviced communities before continuing to a new community the following day. The objective is to minimize the sum of the travel expenses, technician work-hour costs, penalties for community waiting time, and overtime costs based on delays in returning to the depot. To address this problem, they propose a metaheuristic based on the invasive weed optimization algorithm. The authors in~\cite{yahiaoui2023enhanced} address the TRSP with the constraints related to the availability of spare parts and necessary tools, aiming to minimize the overall operational cost. To solve the problem, they develop an enhanced iterated local search algorithm that integrates local search operators, removal heuristics, a best insertion strategy, and an intensification/diversification mechanism based on an elite set of solutions. The dynamic variant of the TRSP is addressed by~\cite{nielsen2023tactical}, where some tasks are known in advance while others arrive dynamically. The authors adopt a tactical planning approach that schedules the known tasks in a way that minimizes overall driving distance, while preserving flexibility to ensure short response times for dynamically arriving tasks. 

As discussed above, existing variants of the TRSP have been formulated with a single objective function. Since our problem involves a lexicographic objective function and is addressed using CG, we next review studies that apply CG to lexicographic optimization in different application domains.

In~\cite{GamachePBS1998, achour2007}, the authors address the preferential bidding system problem, which involves building crew schedules while maximizing their scores in strict order of seniority. Each crew member is assigned an objective function representing their score, which is then optimized lexicographically based on seniority. The authors explain that a weighting approach---where the problem is transformed into a single-objective function using weights---cannot be applied because the weights must increase exponentially with the number of crew members to preserve the lexicographic order. This results in excessively large weights that cannot be handled by any machine for realistic instances. Consequently, they conclude that a sequential approach is necessary. In this approach, the problem is first solved with the objective function of the most senior crew member, then it is solved with the objective function of the second most senior, without reducing the objective of the most senior crew member, and so on. Since solving the problem with integrality constraints at each iteration is computationally expensive, certain integrality conditions are relaxed, which may lead to time-consuming backtracking to ensure feasibility. In~\cite{katrin2022}, the authors also apply a sequential approach to a bin-packing problem with five objective functions, optimized in lexicographic order. They develop a branch-and-price algorithm for each single-objective stage. In~\cite{tellache2025variant}, the authors proposed a weighted approach for a variant of the minimum path cover problem with two objectives. The primary objective is to maximize the number of covered nodes, while the secondary objective is to cover this maximum number of nodes using a minimum number of node-disjoint feasible paths. The problem is solved using a branch-and-cut approach. More recently,~\cite{tellache2024linear} introduces a lexicographic CG approach that considers all objectives simultaneously when generating new columns. This approach, applied to the preferential bidding system, yields promising results compared to the sequential approach, especially in cases involving backtracking. It also avoids the numerical problems caused by the use of large weights.

\section{Problem Description} \label{sec:problem_def}

We denote by $\mathcal{I}$ the set of technical interventions that may be performed on a given work day. Each intervention $i$ has a duration $d_i$ and a time window $[s_i,e_i]$, where $s_i$ and $e_i$ represent the earliest start and latest completion times, respectively. An intervention is classified as \textit{long} if $d_i \geq 120$ minutes and \textit{short} otherwise. We introduce the binary indicator $\theta_i$, which equals 1 if intervention $i$ is long and $0$ otherwise.

To perform the interventions, technicians are grouped into teams, which may consist of a single technician, and each team is assigned a vehicle. For this reason,we use \textit{vehicle} and \textit{team} interchangeably and denote the set of vehicles\slash teams by $\mathcal{V}$. 
The subset $\mathcal{I}_v \subseteq \mathcal{I}$ denotes the set of interventions that vehicle $v\in \mathcal{V}$ can perform, based on its technicians’ skills.

Each vehicle $v$ follows a route that begins and ends at its designated depot, denoted by $\sigma_v$. However, not all vehicles need to operate on a given day; some teams may not work and their vehicles may remain at the depot. The set of all depots is defined as $\mathcal{D} = \bigcup_{v\in \mathcal{V}} \{\sigma_v\}$. Consequently, the set of nodes (i.e., physical locations) in the network is given by $\mathcal{N} = \mathcal{I} \cup \mathcal{D}$, while the nodes that vehicle $v$ can visit are $\mathcal{N}_v = \mathcal{I}_v \cup \{\sigma_v\} \subseteq \mathcal{N}$. The network is fully connected with bidirectional arcs between all nodes $i,j \in \mathcal{N}$; $\ell_{ij}$ denotes the distance between them and $t_{ij}$ the travel time, which satisfies the triangle inequality.

The work day is divided into two halves---morning and afternoon---with a lunch break in between. Since all technicians take lunch at the same time and interventions cannot be scheduled during this period, the endpoint of the morning ($\text{MD}$) coincides exactly with the start of the afternoon. Interventions are non-preemptible, except for long interventions, which may span the lunch break (i.e., the lunch break is the only allowed interruption). Short interventions must be completed within a single half-day. Vehicles must return to their designated depots before the end of the afternoon ($\text{ED}$). 

	As mentioned above, a work day encompasses both work and travel time. While the total duration of each day is globally bounded, additional constraints may be imposed on a set $\mathcal{K}$ of subintervals into which the work day is divided, in order to further limit the workload of a team within each subinterval. To model these constraints in a generic way, we treat subintervals as resources: $q_{ik}\geq 0$ represents the amount of resource $k$ (i.e., time within subinterval $k$) consumed by intervention $i$, and $m_{vk}$ denotes the capacity of vehicle $v$ on resource $k$. Note that $q_{ik}>0$ only when the time window of intervention $i$ lies within subinterval $k$. These additional resource constraints were introduced by EDF to allow buffer adjustments within each interval---by reducing the capacities $m_{vk}$---since the model does not explicitly account for real-world uncertainties in travel or intervention durations.

	The goal is to determine feasible vehicle routes that optimize, in strict lexicographic order:
	(1) the total duration of interventions performed during the day (to be maximized), and (2) vehicle-related costs, including personnel expenses and distance-based travel costs (to be minimized). The total personnel cost associated with vehicle $v$ is denoted by $g_v$, and the cost per unit of distance traveled (per kilometer) is represented by $\delta$. A summary of the introduced notation is provided in Table~\ref{tab:notations}.
	
	\begin{table}[ht]
		
		\caption{Summary of the notation introduced in the problem definition.}
		\label{tab:notations}
		\begin{center}
			\begin{tabular}{ | p{1.7cm} | p{13.5cm}| } 
				\hline
				\textbf{Notation} & \textbf{Description} \\
				\hline
				$\mathcal{I}$ & Set of interventions \\
				$\mathcal{V}$ & Set of vehicles (teams of technicians) \\
				$\mathcal{I}_v \subseteq \mathcal{I}$ & Set of interventions that can be done by vehicle $v \in \mathcal{V}$ \\
				$\mathcal{V}_i \subseteq \mathcal{V}$ & Set of vehicles that can do intervention $i \in \mathcal{I}$ \\
				$\mathcal{D}$ & Set of vehicle depots \\
				$\mathcal{N}$ & Set of nodes ($\mathcal{N} = \mathcal{I} \cup \mathcal{D}$)\\
				$\sigma_v$ & Depot associated with vehicle $v \in \mathcal{V}$ \\
				$\mathcal{N}_v \subseteq \mathcal{N}$ & Set of nodes that can be visited by vehicle $v \in \mathcal{V}$ ($\mathcal{N}_v = \mathcal{I}_v \cup \{\sigma_v\} $)\\
				$\mathcal{K}$ & Set of intervals (not necessarily disjoint) of the work day\\
				$d_i$ & Duration (in minutes) of intervention $i \in \mathcal{I}$ \\
				$s_i$ & Earliest start time of intervention $i \in \mathcal{I}$ \\
				$e_i$ & Latest completion time of intervention $i \in \mathcal{I}$ \\
				$\theta_i$ & Indicator that is 1 if intervention $i \in \mathcal{I}$ is long (i.e., $d_i \geq 120$) and 0 otherwise \\
				$\ell_{ij}$ & Distance between $i \in \mathcal{N}$ and $j \in \mathcal{N}$ (in km)\\ 
				$t_{ij}$ & Travel time between $i \in \mathcal{N}$ and $j \in \mathcal{N}$ (in minutes)\\
				$\text{MD}$ & End of the morning half of the work day \\
				$\text{ED}$ & End of the afternoon half of the work day \\  
				$q_{ik}$ & Time consumption of intervention $i \in \mathcal{I}$ on resource $k \in \mathcal{K}$\\
				$m_{vk}$ & Capacity of vehicle $v \in \mathcal{V}$ and  resource $k \in \mathcal{K}$ \\
				$g_v$ & Daily technician cost associated with vehicle $v \in \mathcal{V}$ (in km)\\
				$\delta$ & Cost per unit distance (\euro \slash km)\\
				\hline
			\end{tabular}
		\end{center}
	\end{table}
	
\section{Mathematical formulation} \label{sec:math_formulation}
\label{sec:math_form}
In this section, we introduce two mathematical formulations of our problem: a compact formulation (Section~\ref{subsec:compactForm}) and an extended formulation based on Dantzig–Wolfe decomposition (Section~\ref{sec:decomposition}).
\subsection{Compact Formulation} \label{subsec:compactForm} 
To formulate the problem as a compact MILP model, we define the following binary variables: \begin{align*}
	x_{ijv} & = \left\{ \begin{array}{rl}
		1& \text{if vehicle }v\text{ goes from node }i \text{ to node }j,\\
		0& \text{otherwise,}
	\end{array}\right . \forall  i,j \in \mathcal{N}_v, v\in \mathcal{V}, i \neq j, \\
	z_i & = \left\{ \begin{array}{rl}
		1& \text{if intervention }i\text{ starts in the afternoon,} \\
		0& \text{otherwise,}
	\end{array}\right . \forall i \in \mathcal{I} : \theta_i = 0, 
\end{align*}
and the continuous variable: \begin{align*}
	u_i \in \mathbb{R}_{\geq 0}:  & \text{ start time of intervention } i, & \forall i \in \mathcal{I}.
\end{align*} Note that $z_i$ is defined only for short interventions (i.e., $\theta_i = 0$), that must start and finish within the same half-day.

The objective function consists of two components that must be considered in lexicographical order. The first objective function, to be maximized, is the total duration of interventions assigned across all vehicles: \begin{equation}
	f_c^1 =  \sum_{v \in \mathcal{V}}\sum_{i \in \mathcal{I}_v} \sum_{j \in \mathcal{N}_v} d_i \cdot x_{ijv}.
\end{equation} The second objective function minimizes the operational cost, which includes both the travel costs between locations and the fixed costs of using vehicles: \begin{equation}
	f_c^2 =  \sum_{v \in \mathcal{V}}\sum_{i \in \mathcal{N}_v} \sum_{j \in \mathcal{N}_v} \ell_{ij} \cdot \delta \cdot x_{ijv} +  \sum_{v \in \mathcal{V}}\sum_{i \in \mathcal{N}_v}  g_v \cdot x_{\sigma_{v}iv}. 
\end{equation}

The compact MILP formulation can be written as follows. \begin{subequations} \label{form:mip_price_conf}
	\begin{align} 
		&&& \text{lexmax} \begin{pmatrix} f_c^1 \\ -f_c^2 \end{pmatrix} \label{eq:CObj} \\
		\text{s.t. } &&&\sum_{v \in \mathcal{V}_i}\sum_{j \in \mathcal{N}_v \setminus \{i\}}
		x_{ijv}\ \le\ 1 && \forall i\in \mathcal{I}  &&  \label{eq:CR1} \\
		&&& \sum_{j \in \mathcal{N}_v, j \neq i} (x_{ijv} - x_{jiv}) = 0 && \forall v \in \mathcal{V}, i\in \mathcal{N}_v  &&  \label{eq:CR2}\\
		&&& \sum_{j \in \mathcal{N}_v} x_{\sigma_v j v} \leq 1 && \forall v \in \mathcal{V} && \label{eq:CR3} \\
		&&& u_i-e_i+(e_i + d_i +t_{ij})\cdot \sum_{v \in \mathcal{V}_i\cap \mathcal{V}_j} x_{ijv} \leq  u_j && \forall i,j\in \mathcal{I}, i \neq j && \label{eq:CR4} \\
		&&& t_{\sigma_v i} \cdot x_{\sigma_v i v} \leq u_i  && \forall v \in \mathcal{D}, i\in \mathcal{I}_v && \label{eq:CR5} \\
		&&& u_i + d_i + t_{i \sigma_v} \cdot x_{i \sigma_v v} \leq \text{ED}  && \forall  v \in \mathcal{V} , i\in \mathcal{I}_v && \label{eq:CR6} \\
		&&& s_i \leq u_i && \forall i\in \mathcal{I} && \label{eq:CR7} \\
		&&& u_i + d_i \leq e_i && \forall i\in \mathcal{I} && \label{eq:CR8} \\
		&&& u_i + d_i \leq \text{MD} + (\text{ED}-\text{MD}) \cdot z_i && \forall i\in \mathcal{I} : \theta_i = 0 && \label{eq:CR9} \\
		&&& u_i \geq \text{MD} \cdot z_i && \forall i\in \mathcal{I} : \theta_i = 0 && \label{eq:CR10} \\
		&&& \sum_{i\in \mathcal{I}_v} \sum_{\substack{j \in \mathcal{N}_v \\ j \neq i}} q_{ik} \cdot x_{ijv} \leq m_{vk} && \forall  v \in \mathcal{V}, k\in \mathcal{K} && \label{eq:CR11}\\
		&&& x_{ijv} \in \{0, 1\} && \forall v \in \mathcal{V}, i,j \in \mathcal{N}_v,\, i \neq j  && \label{eq:domain1} \\
		&&& z_i \in \{0, 1\} && \forall i \in \mathcal{I} : \theta_i = 0 && \label{eq:domain2} \\
		&&& u_i \in \mathbb{R}_{\geq 0} && \forall i \in \mathcal{I} && \label{eq:domain3}
	\end{align}
\end{subequations}

Objective~\eqref{eq:CObj} implements the lexicographic optimization: the model first maximizes the total duration of completed interventions and, among all optimal solutions for this criterion, selects the one with minimum operational cost. Constraints~\eqref{eq:CR1} ensure that each intervention is performed at most once.
Constraints~\eqref{eq:CR2} enforce route consistency for each vehicle by preserving the flow entering and leaving every node. 
Constraints~\eqref{eq:CR3} require each vehicle to depart from its associated depot at most once. 
Constraints~\eqref{eq:CR4} impose temporal feasibility: if intervention~$j$ follows intervention~$i$ on a vehicle’s route, then the start time of intervention~$j$ must be no earlier than the end of intervention~$i$ plus the corresponding travel time.
Constraints~\eqref{eq:CR5} define the earliest feasible start time of the first intervention assigned to a vehicle. Note that constraints~\eqref{eq:CR4} and~\eqref{eq:CR5} prevent subtours, as the start times of interventions are propagated through the precedence relations between them. Combined with the flow constraints~\eqref{eq:CR2}, they further ensure that each vehicle must pass through its designated depot. In particular, because the depot does not appear in the temporal linking constraints~\eqref{eq:CR4}, and every visited node must satisfy the flow-balance conditions in~\eqref{eq:CR2}, the first and last interventions assigned to a vehicle must be connected to the depot. Constraints~\eqref{eq:CR6} limit the total working and travel time to the duration of the work day (i.e., until $\text{ED}$). Constraints~\eqref{eq:CR7}~and~\eqref{eq:CR8} ensure compliance with the time window associated with each intervention. Constraints~\eqref{eq:CR9}~and~\eqref{eq:CR10} enforce that $z_i = 1$ if and only if a short intervention $i$ begins in the afternoon.
Finally, constraints~\eqref{eq:CR11} require that, for each vehicle and each interval, the total amount of consumed resource does not exceed the available capacity.

\subsection{Extended formulation}\label{sec:decomposition}

We also propose an extended formulation based on the well-known set-packing formulation often used in routing problems. In this formulation, variables representing trips between pairs of nodes are replaced by variables representing the entire routes. In our case, these route variables additionally encode the timing of interventions and an indicator for whether a given intervention is scheduled in the afternoon. These features are captured by the variables $u_i$ and $z_i$ from the compact MILP formulation, respectively.

Formally, let $\mathcal{R}_v$ denote the set of all feasible routes that vehicle $v$ can perform in a single day, taking into account location constraints, time windows, the lunch break and resource capacities. Each route $r \in \mathcal{R}_v$ thus corresponds to a feasible assignment of variables $x$, $u$, and $z$ that satisfy constraints (\ref{eq:CR2})--(\ref{eq:CR11}). In the set-packing formulation, the original variables $x$, $u$, and $z$ are replaced by a new set of binary variables  $\{\lambda_{rv} \mid r \in \mathcal{R}_v, v \in \mathcal{V} \}$, where $\lambda_{rv} = 1$ if and only if route $r$ is selected for vehicle $v$ in the final solution. 

The variables $\lambda$ capture key information associated with the variables $x$, $u$, and $z$ from the compact formulation. For a given route $r$ and vehicle $v$, we precompute the following data:
\begin{itemize}
	\item $a_{ir} = 1$ if intervention $i$ is included in route $r$, and $0$ otherwise;
	\item $D_r$: the total duration, in minutes, of all interventions included in route $r$;
	\item $c_r$: the total cost associated with route $r$, including both travel and technician costs.
\end{itemize}

The objective functions $f_c^1$ and $f_c^2$ can then be expressed in terms of the $\lambda$ variables: \begin{subequations}
	\begin{align}
		f_e^1 & =  \sum_{v \in V} \sum_{r \in \mathcal{R}_v} D_r \cdot \lambda_{rv},\\
		f_e^2 & =  \sum_{v \in V} \sum_{r \in \mathcal{R}_v} c_r \cdot \lambda_{rv}.        
	\end{align}
\end{subequations} Using this reformulation, the model (\ref{eq:CObj})--(\ref{eq:domain3}) can be expressed as: 
\begin{subequations}
	\begin{align} 
		& \text{lexmax} \begin{pmatrix} f_e^1\\ -f_e^2 \end{pmatrix}  \label{DW:Obj} \\
		\text{s.t.} \nonumber \\   
		& \sum_{v \in \mathcal{V}} \sum_{r \in \mathcal{R}_v}   a_{ir} \lambda_{rv}  \leq 1 && \forall i \in \mathcal{I} \label{DW:cover}\\
		&  \sum_{r \in \mathcal{R}_v} \lambda_{rv} \leq 1 &&  \forall v \in \mathcal{V} \label{DW:onerouteperv}\\ 		
		& \lambda_{rv} \in \{0,1\} && \forall v \in \mathcal{V}, \quad r \in \mathcal{R}_v \label{DW:routevar}
	\end{align}
\end{subequations}

Formulation~(\ref{DW:Obj})--(\ref{DW:routevar}) involves an exponential number of route variables. Nevertheless, such formulations generally provide stronger linear relaxations that can be solved efficiently by dynamically generating route variables using a CG-based procedure, which is described in Section~\ref{sec:ColGen}.

\section{Single-objective reformulations of the lexicographic TRSP} \label{sec:solution_methods}


Classical approaches for solving multi-objective optimization problems typically transform them into a sequence of single-objective problems. A common technique is the weighted-sum method \cite{marler2010weighted}, which aggregates multiple objectives into a single scalar objective. Another common approach is the $\epsilon$-constrained method \cite{haimes1971bicriterion}, which optimizes one objective while enforcing bounds on the others (i.e., treating them as constraints).

Lexicographic optimization, in contrast, treats objectives according to a strict hierarchical order. A lexicographically optimal solution dominates all others with respect to this predefined hierarchy. When classical techniques such as the weighted-sum method are applied in this context, the goal is not to explore trade-offs but to determine weights that enforce the hierarchy. In this section, we adapt two approaches to handle the lexicographic objective function of our TRSP. In \Cref{subsec:weighted}, we detail how the weighted-sum method can be modified to respect lexicographic priorities through appropriate weight selection. In \Cref{sub:sequential}, we present a sequential approach inspired by the $\epsilon$-constrained method, which optimizes objectives in priority order while preserving the optimal value of higher-priority objectives.

	\subsection{Weighted approach}  \label{subsec:weighted}
	
	
	In this approach, the objective function combines both components into a single term using a big-$M$ constant as a weighting factor: \begin{align}
		\text{max }  &&& M \cdot f^1 - f^2. && \label{eq:WA}
	\end{align} The value of $M$  must ensure that maximizing the total amount of time dedicated to interventions is always prioritized over minimizing vehicle costs. As will be discussed in \Cref{sub:sequential}, $M$ should be chosen as small as possible while still enforcing the required strict hierarchy between objectives. The value we employ is \begin{align}\label{eq:M}
		M & = \left[ |\mathcal{V}| \cdot \left( \text{ED} - \min_{i \in \mathcal{I}} d_i \right) \cdot \max_{i,j \in \mathcal{N}, i \neq j}\frac{\ell_{ij}}{t_{ij}} \cdot \delta + \left( \sum_{v \in \mathcal{V}} g_v - \min_{v \in \mathcal{V}} g_v \right) \right] \cdot \frac{1}{\gcd \{ d_i, \forall i \in \mathcal{I} \} }  +1.
	\end{align}
	
	For convenience, let $I=|\mathcal{I}|$ and $V=|\mathcal{V}|$, and denote by $M_0$ the expression inside the square brackets in \eqref{eq:M}. We now argue that $M_0$ provides an upper bound on the possible increase in vehicle costs associated with a one-unit increase in total intervention time. 
	
	If the increase is from zero to one unit of intervention time, the bound is trivially valid, as only one vehicle is required. Otherwise, we may assume that at least one intervention is being performed and one vehicle is already in use in the incumbent solution. The first term of $M_0$ corresponds to the worst-case additional operational cost: all $V$ vehicles operate for the entire day while performing at least one intervention, each traveling along the most expensive street segment per time unit. The second term represents the maximum possible increment in personnel costs. Because intervention durations can only increase in multiples of $\gcd \{ d_i, \forall i \in \mathcal{I}\}$, and the objective multiplies the total intervention time by $M$, we may safely divide $M_0$ by this $\gcd$, as shown in equation~(\ref{eq:M}). Finally, we add 1 to this value to ensure strict lexicographic order between objectives. 
	
	To illustrate that the upper bound $M_0$ can in principle be attained, consider the following worst-case scenario. Let $I = V$, and assume that every intervention has duration $d_i = 1$ for all $i \in \mathcal{I}$, and that all vehicles have identical personnel cost $g_v = g$. We assume that vehicle $v = 1$ is capable of performing any intervention, whereas each vehicle $v \in \{2,\ldots,V\}$ can perform only intervention $i = v$. The geographical layout is such that $t_{\sigma_i i} = t_{i \sigma_i} = \frac{\mathrm{ED} - 1}{2}$ for all $i \in \mathcal{I}$, while $t_{\sigma_1 i} = t_{i \sigma_1} = 0$ for all $i \in \{2,\ldots,I\}$; moreover, all ratios $\frac{\ell_{ij}}{t_{ij}}$ are identical. Consider solution A, in which vehicle $1$ executes all interventions except $i = 1$, yielding a total intervention duration of $V - 1$ and a vehicle cost of $g$. In contrast, solution B assigns each vehicle $v$ to intervention $i = v$, achieving a total intervention duration of $V$ but incurring a total vehicle cost of $M_0 + g$. Hence, switching from solution A to solution B increases the completed intervention duration by exactly one unit, while increasing vehicle-related costs by precisely $M_0$.

	\subsubsection{Compact Formulation} \label{subsec:compact}
	
	In the compact formulation, the weighted objective function~\eqref{eq:WA} translates into the expression given in~\eqref{eq:OF1}. The first term rewards executed interventions in proportion to their durations, multiplied by the constant \(M\), which enforces the strict priority of maximizing completed intervention time. The second term penalizes vehicle usage, capturing both travel-related operational costs and personnel costs. The constant \( M \) is computed from the input data and ensures that intervention coverage is prioritized over cost minimization. \begin{align}
		\text{max }  &&&  \sum_{v \in \mathcal{V}}\sum_{i \in \mathcal{I}_v} \sum_{\substack{j \in \mathcal{N}_v \\ j \neq i}} M \cdot d_i \cdot x_{ijv} 
		- \left( \sum_{v \in \mathcal{V}} \sum_{i \in \mathcal{N}_v} \sum_{\substack{j \in \mathcal{N}_v \\ j \neq i}} \ell_{ij} \cdot \delta \cdot x_{ijv} 
		+ \sum_{v \in \mathcal{V}}\sum_{i \in \mathcal{N}_v}  x_{\sigma_{v}iv} \cdot g_v \right) \label{eq:OF1} \\
		\text{s.t. } &&& \eqref{eq:CR1}-\eqref{eq:domain3} \nonumber
	\end{align}

	\subsubsection{Extended Formulation} \label{subsub:extended_form_weighted}
	
	In the extended formulation, the weighted objective function~\eqref{eq:WA} is expressed in terms of route variables, as shown in~\eqref{SP:obj}.  Each selected route \( r \in \mathcal{R}_v \) contributes \( D_r \) units of intervention time and incurs a routing cost \( c_r \). The term \( M \cdot D_r \) enforces the lexicographic priority of maximizing completed intervention time, while the second term aggregates all routing costs. As in the compact formulation, the constant \( M \) guarantees that intervention coverage dominates cost minimization. \begin{align} 
		\text{max }  &&& \sum_{v \in \mathcal{V}} \sum_{r \in \mathcal{R}_v} M \cdot D_{r} \cdot \lambda_{rv} - \sum_{v \in \mathcal{V}} \sum_{r \in \mathcal{R}_v} c_{r} \cdot \lambda_{rv} \label{SP:obj} \\
		\text{s.t. } &&& \eqref{DW:cover}-\eqref{DW:routevar} \nonumber
	\end{align}
\subsection{Sequential Approach}
\label{sub:sequential}

A weighted formulation allows us to express the problem as a MILP model with a single objective function, but it requires the use of a data-dependent weighting coefficient \(M\). Since \(M\) must enforce strict lexicographic priority, it often becomes very large in practice, which may severely affect numerical performance. As is well known, such large coefficients can cause numerical instability in MILP solvers, leading to poor convergence, rounding errors, and loss of precision (see, e.g., \cite{cococcioni_big-m_2021}). These difficulties become even more pronounced in our extended formulation, where the problem is solved via a CG scheme. Similar issues have been documented in the literature. For example, \cite{baldacci2011exact} studies a pickup-and-delivery problem with time windows in which a very large fixed vehicle cost \(W\) (ensuring that minimizing the fleet size takes precedence over routing cost) leads to numerical instability and slow convergence within the CG process.

To overcome these limitations, we propose a sequential optimization approach that avoids the use of big-$M$ coefficients altogether. In this method, objectives are optimized one at a time according to their predefined priority. At each stage, the current objective is optimized while preserving the optimal values already achieved for all higher-priority objectives.

\subsubsection{Compact Formulation}

In the first phase, the total duration of performed interventions is maximized: \begin{align}
	\text{max} \quad & f_c^1 = \sum_{v \in \mathcal{V}}  \sum_{i \in \mathcal{I}_v} \sum_{\substack{j \in \mathcal{N}_v \\ j \neq i}} d_i \cdot x_{ijv} \label{eq:OF2a} \\
	\text{s.t.} \quad & \eqref{eq:CR1} - \eqref{eq:domain3} \nonumber
\end{align} 

Let \( d^{*}_c \) denote the optimal value of \( f_c^1 \). In the second phase, the secondary objective minimizes the total operational cost while enforcing the primary objective’s optimal value: \begin{subequations}
	\begin{align}
		\text{min} \quad & f_c^2 = \sum_{v \in \mathcal{V}} \sum_{i \in \mathcal{N}_v} \sum_{j \in \mathcal{N}_v} \ell_{ij} \cdot \delta \cdot x_{ijv} 
		+ \sum_{v \in \mathcal{V}}\sum_{i \in \mathcal{N}_v}  g_v \cdot x_{\sigma_{v}iv} \label{eq:OF2b} \\
		\text{s.t.} \quad & \eqref{eq:CR1} - \eqref{eq:domain3} \nonumber \\
		& \sum_{v \in \mathcal{V}}\sum_{i \in \mathcal{I}_v} \sum_{\substack{j \in \mathcal{N}_v \\ j \neq i}}  d_i \cdot x_{ijv} \geq d_c^* \label{eq:lexConstraint}
	\end{align}
\end{subequations} Constraint~\eqref{eq:lexConstraint} ensures that the total duration of completed interventions remains at least $d_c^*$, thereby preserving the optimal value of the primary objective.

\subsubsection{Extended Formulation}
\label{subsub:extended_form_sequentiel}

Similarly, in the first phase, the total duration of completed interventions is maximized using route-based variables: \begin{align}
	\text{max} \quad & f_e^1 = \sum_{v \in \mathcal{V}} \sum_{r \in \mathcal{R}_v} D_r \cdot \lambda_{rv} \label{SF1:obj} \\
	\text{s.t.} \quad & \eqref{DW:cover} - \eqref{DW:routevar} \nonumber
\end{align}

Let \( d^{*}_e \) denote the optimal value of \( f_e^1 \). In the second phase, the secondary objective minimizes the total cost of the selected routes while maintaining the primary objective’s optimal value: \begin{subequations}
	\begin{align}
		\text{min} \quad & f_e^2 = \sum_{v \in \mathcal{V}} \sum_{r \in \mathcal{R}_v} c_r \cdot \lambda_{rv} \label{SF2:obj} \\
		\text{s.t.} \quad & \eqref{DW:cover} - \eqref{DW:routevar} \nonumber \\
		& \sum_{v \in \mathcal{V}} \sum_{r \in \mathcal{R}_v} D_r \cdot \lambda_{rv} \geq d^{*}_e  \label{SF2:duration}
	\end{align}
\end{subequations} 

Constraint~\eqref{SF2:duration} plays the same role as its compact counterpart: it ensures that any reduction in total cost does not compromise the duration of interventions achieved in the first phase.

\section{Column-generation-based algorithm} \label{sec:ColGen}

In \Cref{sec:solution_methods}, we introduced single-objective reformulations of the considered lexicographic TRSP variant, using both weighted and sequential approaches. The resulting compact formulations can be solved directly with (commercial) MILP solvers, whereas the extended formulations are typically addressed through methods based on CG~\cite{cortes2014branch}, which we describe in this section.

Let us denote the extended MILP formulation by P. CG solves its master problem (MP), i.e., the linear relaxation of~P, without explicitly enumerating all variables \( \lambda_{rv} \). The algorithm starts with a restricted subset $\mathcal{R}' \subseteq \cup_{v \in \mathcal{V}} \mathcal{R}_v$ of variables (columns), which ---by slight abuse of language---we simply refer to as routes (see \cref{sec:decomposition}). The restricted master problem MP($\mathcal{R}'$) is then solved. Using the optimal dual solution, the algorithm searches for routes not yet in $\mathcal{R}'$ that have positive reduced cost (since we maximize). If at least one such route is found, it is added to \( \mathcal{R}' \) and the process is repeated; otherwise, MP has been solved to optimality.

The computation of reduced costs relies on the optimal dual solution of MP($\mathcal{R}'$) and on the contribution of the route to the objective function. Identifying new routes with positive reduced cost constitutes the pricing subproblem (PS). Each PS is solved independently for a vehicle $v$ to identify routes $r$ with positive reduced cost $\bar{c}_r$. We describe the PS and our solution method in more detail in Section \ref{subsec:PSModelingSolving}. 

Once MP has been solved to optimality, we solve P($\mathcal{R}'$), i.e., the restriction of P to the current set of routes \( \mathcal{R}' \). Let \( u \) be the optimal objective value of MP and \( \ell \) the optimal objective value of P($\mathcal{R}'$). It has been shown in~\cite{nemhauser1988,crowder1983}, and was already used by~\cite{dantzig1954solution} to solve to optimality a 49-city traveling salesman problem, that any column with reduced cost less than $\ell - u$ cannot belong to an optimal solution. Therefore, it is sufficient to generate only columns with reduced cost $\bar{c}_r \ge \ell - u$ to ensure optimality. If we add those columns to $\mathcal{R}'$, then solving $P(\mathcal{R}')$ is equivalent to solving P, i.e., the full extended formulation. The resulting procedure is summarized in Algorithm~\ref{Exact_algo}.

\begin{algorithm}[H]
	\caption{Exact CG-based algorithm} \label{Exact_algo}
	\begin{algorithmic}[1]
		\State \textbf{Initialize:} $\mathcal{R}' \gets$ initial set of routes; \label{Exact_algo_Init}
		\Repeat
		\State Solve MP($\mathcal{R}'$); \Comment{Let the optimal value be $u$}
		\State For each vehicle \( v \), solve the PS to generate routes $r$ with \( \bar{c}_r > 0 \) and add them to \( \mathcal{R}' \);
		\Until{No feasible route with \( \bar{c}_r > 0 \) exists}\label{Exact_algo_end_CG_loop}
		\State Solve P($\mathcal{R}'$); \Comment{Let the solution value be $\ell$}
		\State For each vehicle \( v \), solve the PS to generate all feasible routes with \( \bar{c}_r \geq \ell - u \) and add them to \( \mathcal{R}' \);
		\State Solve P($\mathcal{R}'$);
		\State \Return Optimal solution of P (extended formulation).
	\end{algorithmic}
\end{algorithm}

\subsection{Adaptations to the lexicographic extended formulation}

In this section, we show how Algorithm~\ref{Exact_algo} is
adapted to the two single–objective formulations introduced in
\cref{sec:solution_methods}. Specifically, we detail how the reduced costs and PS change under the weighted and the sequential approaches.


\subsubsection{Weighted approach} \label{subsub:weigthed_approach_with_lex}

In the weighted extended formulation, the two lexicographic objectives are linearly combined in a single expression. Hence, the problem can be solved in one run of Algorithm~\ref{Exact_algo}. The reduced cost for a route $r\in\mathcal{R}_v$:
\[
\bar{c}_r = M\,D_r - c_r - \mu_v - \sum_{i\in \mathcal{I}} a_{ir}\,\nu_i ,
\]
where $D_r$ and $c_r$ are respectively the duration and cost of route $r$, $a_{ir}$ is a binary indicator equal to 1 if intervention \( i \) is covered by route \( r \), and 0 otherwise, and $(\mu_v,\nu_i)$ are respectively the optimal values of the dual variables associated with constraints \eqref{DW:onerouteperv}--\eqref{DW:cover}. This formulation ensures that route generation remains consistent with the weighted prioritization of objectives.



\subsubsection{Sequential approach} \label{subsub:sequentiel_with_lex}

In the sequential extended formulation, the two objectives are optimized in two successive executions of Algorithm~\ref{Exact_algo}. The first phase maximizes the total intervention duration, using the reduced cost
\[
\bar{c}_r^1 = D_r - \mu_v^1 - \sum_{i \in \mathcal{I}} a_{ir}\nu_i^1 ,
\]
where $\mu_v^1$ and $\nu_i^1$ are the optimal values of the dual variables of the first-phase MP($\mathcal{R}'$). Routes with $\bar{c}_r^1>0$ are iteratively generated until the first-phase MP($\mathcal{R}'$) is solved to optimality.

The second phase minimizes the total routing cost while enforcing the previously attained intervention duration. The PS uses the reduced cost
\[
\bar{c}_r^2 = c_r - \mu_v^2 - \sum_{i \in \mathcal{I}} a_{ir}\nu_i^2 - D_r\,\rho ,
\]
where $\mu_v^2$ and $\nu_i^2$ are the optimal values of the dual variables of the second-phase MP($\mathcal{R}'$) and $\rho$ is the optimal value of the dual variable associated with the minimum-duration constraint~\eqref{SF2:duration}.

\subsection{Modeling and solving the pricing subproblems} \label{subsec:PSModelingSolving}



We model the PS as an ESPPRC; see \cref{sec:LR} for references on this problem. For each vehicle~\( v \), we consider a complete directed graph \( G_v = (\mathcal{N}_v, \mathcal{A}_v) \), where \( \mathcal{A}_v \) contains all arcs whose endpoints lie in \( \mathcal{N}_v \). Each arc \( (i,j) \in \mathcal{A}_v \) is associated with a travel time \( t_{ij} \) and with a cost chosen so that the total cost of any path from \( \sigma_v \) back to \( \sigma_v \) coincides with the reduced cost of the corresponding variable \( \lambda_{rv} \) (multiplied by -1 in the case of maximization). The ESPPRC then consists in finding, in \( G_v \), a maximum-reduced-cost path that satisfies all operational constraints: compliance with time windows, respect for vehicle-capacity limits along the route, and the requirement that each intervention node in \( I_v \subseteq \mathcal{N}_v \) is visited at most once.

We solve the ESPPRC using a labeling algorithm based on dynamic programming. This method incrementally constructs feasible paths and propagates associated labels, while relying on dominance rules to discard inferior labels and keep the search tractable. The full details of the algorithm are presented in the following subsections. For clarity, we present the algorithm for the maximization case, but the adaptation to minimization is straightforward.

\subsubsection{Labeling}
We use a label $L_{v,i}$ to represent a partial path in the graph $G_v$ that terminates at node $i\in \mathcal{N}_v$. Each label stores the following information: $\bar{c}(L_{v,i})$, the reduced cost of the partial path; $t(L_{v,i})$, the arrival time at node $i$; $m_{vk}(L_{v,i})$, the cumulative consumption of resource $k$; $\mathcal{U}(L_{v,i}) \subseteq \mathcal{N}_v$, the set of nodes already included in the path; and $\mathcal{P}(L_{v,i}) \subseteq \mathcal{I}_v$, the set of interventions that remain eligible for extension. The objective of the labeling algorithm is to generate a set of feasible routes with positive reduced costs. 

The pseudocode of the labeling procedure is provided in \cref{alg:labeling}. An initially empty list, \textproc{new\_routes}, gathers the feasible routes generated during execution. Two additional lists \textproc{unexpanded\_paths} and \textproc{dominating\_paths} are introduced to keep track of the partial paths that have not yet been expanded and those that dominate others (see \cref{sec:domination}), respectively. At initialization, a single path object is created--corresponding to the depot node \( \sigma_v \)--and inserted into \textproc{unexpanded\_paths}. 
This object includes the path identifier, the vehicle identifier ($v$), the last node in the path ($i=\sigma_v$), the reduced cost ($\bar{c}(L_{v,i})$, which is set to the negative sum of the fixed cost associated with vehicle $v$, i.e., $g_v$ and the optimal value of the dual variable $\mu_v$), the arrival time ($t(L_{v,i})$, initialized to 0), the consumed-capacity variables ($\overline{m}_{vk}$, which are initialized to 0), the nodes visited in the path ($\mathcal{U}(L_{v,\sigma_v})$), the associated traveling cost ($g_v$), the cumulative duration of performed interventions (as given by $\sum_{j \in \mathcal{U}(L_{v,i}), j \neq \sigma_v} d_j$), and the potential interventions that can still be performed ($\mathcal{P}(L_{v,i}) = \mathcal{I}_v$). 


The procedure iterates while \textproc{nodes\_to\_visit} is non-empty; this list initially contains only the depot node (so that the vehicle terminates its route there). 
Nodes in \textproc{nodes\_to\_visit} are sorted in increasing order of a cost–benefit ratio: 
the cost corresponds to the earliest completion time of the intervention, whereas the benefit reflects the increase in reduced cost resulting from its inclusion. 
At each iteration, the first node~\( i \) in the list is selected. 
All partial paths ending at node \( i \) in \textproc{unexpanded\_paths} are then sorted in decreasing order of their reduced cost, 
so that more promising paths are expanded first, which allows to reduce the number of paths to be explored. 
Each of these paths is then considered for extension to every reachable node \( j \neq i \). 
If the extension is feasible, the resulting label is added to \textproc{unexpanded\_paths}, provided it is not dominated by an existing label (see Section~\ref{sec:domination}). 
When an extension returns to the depot node \( \sigma_v \), the corresponding complete route is evaluated and, if it has positive reduced cost, added to \textproc{new\_routes}.

The procedure terminates when all partial paths have been explored or when a user-defined number \( \eta \) of new routes has been generated. 
The resulting set of feasible routes with positive reduced cost is returned and added to \( \mathcal{R}' \) for the next iteration of Algorithm~\ref{Exact_algo}.

\begin{algorithm}[ht!]
	\caption{Labeling algorithm for vehicle $v$}
	\label{alg:labeling}
	\begin{algorithmic}[1]
		\Function{Find Positive Reduced Cost Paths}{$v$}
		\State \textbf{Initialize:}
		\Statex \hspace{0.5cm} $\bar{c}(L_{v,\sigma_v}) \gets - (g_v + \mu_v)$
		\Statex \hspace{0.5cm} $\overline{m}_{vk} \gets 0, \quad \forall k$
		\Statex \hspace{0.5cm} $t(L_{v,\sigma_v}) \gets 0$
		\Statex \hspace{0.5cm} \textproc{new\_routes} $\gets \emptyset$
		\Statex \hspace{0.5cm} \textproc{unexpanded\_paths} $\gets \{ \}$, \textproc{dominating\_paths} $\gets \{ \}$
		\Statex \hspace{0.5cm} \textproc{nodes\_to\_visit} $\gets [\sigma_v]$
		
		\Repeat 
		\State $i \gets \textbf{popfirst!}(\textproc{nodes\_to\_visit})$
		\State \textbf{Sort}(\textproc{unexpanded\_paths}$[i]$, decreasing by $\bar{c}$)
		\State \textbf{Sort}(\textproc{nodes\_to\_visit}, increasing by cost/benefit ratio)
		
		\ForAll{$L_{v,i} \in \textproc{unexpanded\_paths}[i]$}
		\State $L_{v,j} \gets$ \textproc{ExtendPath}($L_{v,i}$, $j \in $ \textproc{nodes\_to\_visit})
		\Statex \Comment{See Section~\ref{sec:extension}}
		\If{$L_{v,j} \neq \text{null}$}
		\If{\textproc{IsDominated}($L_{v,j}$, \textproc{dominating\_paths}) = \textbf{false}} 
		\Statex \Comment{See Section~\ref{sec:domination}}
		\If{$j \neq \sigma_v$}
		\State \textproc{AddToDict}(\textproc{unexpanded\_paths}, $\text{next\_node}$, $L_{v,j}$)
		\State \textproc{AddIfAbsent}(\textproc{nodes\_to\_visit}, $j$)
		\Else
		\If{$\bar{c}(L_{v,j}) > 0$}
		\State \textproc{new\_routes} $\gets \textproc{new\_routes} \cup \{L_{v,j}\}$
		\If{$|\textproc{new\_routes}| \geq \eta$}
		\State \Return \textproc{new\_routes}
		\EndIf
		\EndIf
		\EndIf
		\EndIf
		\EndIf
		\EndFor
		\State \textproc{unexpanded\_paths}[i] $\gets \emptyset$
		\Until{\textproc{nodes\_to\_visit} $= \emptyset$}
		
		\State \Return \textproc{new\_routes}
		\EndFunction
	\end{algorithmic}
\end{algorithm}

\subsubsection{Extension function} \label{sec:extension}

The extension function is a key component of the labeling algorithm, as it is responsible for the generation of feasible paths. When expanding $L_{v,i}$ to $L_{v,j}$, the procedure considers two cases: $j=\sigma_v$ or $j \in \mathcal{I}_v$. In the first case, the label $L_{v,\sigma_v}$ is considered feasible if the arrival time satisfies $t(L_{v,\sigma_v}) \leq \text{ED}$. In the second case, we introduce the set $\mathcal{O}(L_{v,i})$ of all interventions that are no longer feasible after reaching intervention $j$ due to time- or capacity-related restrictions. A label $L_{v,j}$ arriving at an intervention node \(j \in \mathcal{I}_v\) is feasible if \(t(L_{v,j}) + d_j \leq e_j\) (i.e., the intervention can be completed before its latest allowable end time), \(\overline{m}_{vk}(L_{v,j}) \leq m_{vk}~\forall k \in K\) (i.e., the consumed capacity in each interval does not exceed the available capacity) and \(j \notin \mathcal{U}(L_{v,i})  \). Extended labels that do not meet these conditions are discarded. These extension rules are summarized in \Cref{tab:ext_rules}.

\begin{table}[htbp]
	\centering
	\caption{Extension rules}
	\label{tab:ext_rules}
	
	\vspace{0.75em}
	
	\resizebox{\textwidth}{!}{
		\begin{tabular}{|l|l|}
			\hline
			$j=\boldsymbol{\sigma_v}$ & $\boldsymbol{j \in I_v}$ \\
			\hline
			
			$\mathcal{P}(L_{v,j}) = \emptyset$ 
			&
			$\mathcal{P}(L_{v,j}) = \mathcal{P}(L_{v,i}) \setminus \mathcal{O}(L_{v,i})$ 
			\\
			\hline
			
			$t(L_{v,j}) = t(L_{v,i}) + d_j + t_{i\sigma_v}$ 
			&
			\begin{tabular}[c]{@{}l@{}} 
				$t(L_{v,j}) = \left\{ 
				\begin{array}{ll}
					\max(\text{MD}, s_j), & \text{if } \theta_i = 0 \text{ and } \text{MD} - d_j < t(L_{v,i}) + d_i + t_{ij} \leq \text{MD}, \\[4pt]
					\max(t(L_{v,i}) + d_i + t_{ij}, s_j), & \text{otherwise.}
				\end{array} 
				\right.$
			\end{tabular}
			\\
			\hline
			
			$\overline{m}_{vk}(L_{v,j}) = \overline{m}_{vk}(L_{v,i}), \ \forall k \in K$
			&
			$\overline{m}_{vk}(L_{v,j}) = \overline{m}_{vk}(L_{v,i}) + q_{jk},\ \forall k \in K$
			\\
			\hline
		\end{tabular}
	}
	
\end{table}

	
	
		

\subsubsection{Domination} \label{sec:domination}

Dominance rules are essential for reducing the search space of the labeling algorithm by eliminating partial paths that cannot lead to an optimal or improving complete route. Dominance conditions (indexed by $\beta_s, s = 1, \dots, 4$) are used to compare two labels \(L_{v,j}^1\) and \(L_{v,j}^2\), which represent different paths of a given vehicle $v$ reaching node $j$. A label \(L_{v,j}^1\) is said to dominate \(L_{v,j}^2\) if it is no worse in all considered criteria and strictly better in at least one of them. When a label is dominated, it can be safely discarded because extending it cannot lead to a better feasible route. 

Condition $\beta_1$ requires that the reduced cost associated with $L_{v,j}^1$ must be greater than or equal to the reduced cost associated with $L_{v,j}^2$. A label $L_{v,j}^1$ is considered dominant with respect to condition $\beta_2$  if the set of potential interventions that can still be visited by $L_{v,j}^2$ is a subset of those that can still be visited by $L_{v,j}^1$. This means that label $L_{v,j}^1$ has the potential to visit at least all nodes that $L_{v,j}^2$ can visit, and potentially more. In condition $\beta_3$, label $L_{v,j}^1$ dominates label $L_{v,j}^2$ if $L_{v,j}^1$ has an earlier or the same arrival time at node $j$ with respect to $L_{v,j}^2$. Nevertheless, $L_{v,j}^1$ can still dominate  $L_{v,j}^2$ even with a later arrival time, provided that $L_{v,j}^2$ cannot visit any further nodes. Finally, in terms of consumed capacities, condition $\beta_4$ states that  label $L_{v,j}^1$ dominates  label $L_{v,j}^2$ if it has a smaller or equal consumed capacity for all intervals \(k\). However, if $L_{v,j}^1$ has a higher consumption for an interval, it can still dominate $L_{v,j}^2$ if the additional capacity that $L_{v,j}^1$ could consume for all remaining potential interventions does not exceed its remaining capacity, or if $L_{v,j}^2$ has no remaining potential interventions. The dominance conditions are summarized in \Cref{tab:dominance}.

\begin{table}[htbp]
	\caption{Dominance conditions: $L_{v,j}^1$ dominates $L_{v,j}^2$ if the condition is satisfied}%
	\label{tab:dominance}%
	\begin{center}
		
		\begin{tabular}{|c|ll|}
			\hline
			&  \multicolumn{2}{|l|}{\textbf{Dominance conditions}}\\
			\hline
			$\beta_1$ & $ \bar{c}(L_{v,j}^1)$ & $\geq \bar{c}(L_{v,j}^2)$\\
			\hline
			$\beta_2$ &  $\mathcal{P}(L_{v,j}^2)$ & $\subseteq \mathcal{P}(L_{v,j}^1)$ \\
			\hline
			$\beta_3$ &  $ t(L_{v,j}^1)$ & $\leq t(L_{v,j}^2) $ \textbf{ OR }\\
			&  $ t(L_{v,j}^1)$ &$> t(L_{v,j}^2) \text{ and } \mathcal{P}(L_{v,j}^2) = \emptyset $ \\
			\hline
			$\beta_4$ &  $ \overline{m}_{vk} (L_{v,j}^1)$ &$\leq \overline{m}_{vk} (L_{v,j}^2)$ \textbf{ OR }\\
			&  $ \overline{m}_{vk}(L_{v,j}^1) $&$> \overline{m}_{vk}(L_{v,j}^2) \text{ and } \sum_{i \in \mathcal{P}(L_{v,j}^1)} q_{ik} \leq m_{vk} - \overline{m}_{vk}(L_{v,j}^1) $ \textbf{ OR }\\
			&  $ \overline{m}_{vk}(L_{v,j}^1)$&$> \overline{m}_{vk}(L_{v,j}^2) \text{ and } \sum_{i \in \mathcal{P}(L_{v,j}^2)} q_{ik} = 0 $    \\ 
			\hline
		\end{tabular}%
	\end{center}
\end{table}%

\subsubsection{Relaxation of dominance} \label{subsub:relaxation}

Enforcing all the dominance conditions presented in \Cref{sec:domination} from the beginning might be computationally expensive, negatively affecting overall performance. To enhance computational efficiency, we propose a dominance relaxation strategy that gradually enforces the dominance conditions during the search process. The order in which they are enforced corresponds to their identifier $\beta_s$. This strategy allows us to focus the search on the most promising regions of the solution space, thereby reducing computational time without compromising the quality of the generated routes. 

To implement this approach, we replace lines \ref{Exact_algo_Init} to \ref{Exact_algo_end_CG_loop} in Algorithm~\ref{Exact_algo} by  Algorithm~\ref{alg:dominance}. At each iteration, the MP($\mathcal{R}$) is solved, while the PS generates new routes considering only a subset $\{\beta_1, \ldots, \beta_s\}$ of dominance conditions. If no new route with a positive reduced cost is found (i.e., $\mathcal{R}^{\text{new}} = \emptyset$), an additional dominance condition is activated ($s+1$). The process continues until all dominance rules are applied ($s = 4$).

\begin{algorithm}[H]
	\caption{Dominance relaxation algorithm}
	\label{alg:dominance}
	\begin{algorithmic}[1]
		\State $s \gets 1$ \Comment{Identifier of dominance conditions}
		\State $\mathcal{R} \gets$ Initial set of routes
		\Repeat
		\State Solve MP($\mathcal{R}$)
		\State Solve the PS applying the subset of dominance rules $\{ \beta_1, \ldots, \beta_s \}$
		\State $\mathcal{R}_{\text{new}} \gets$ new routes with positive reduced cost
		\If{$\mathcal{R}^{\text{new}} \neq \emptyset$}
		\State $\mathcal{R} \gets \mathcal{R} \cup \mathcal{R}^{\text{new}}$
		\Else
		\State $s \gets s + 1$ \Comment{Activate one additional dominance condition}
		\EndIf
		\Until{$s > 4$}
	\end{algorithmic}
\end{algorithm}

\section{Computational results} \label{sec:comp_result}
This section reports computational experiments evaluating the performance of the solution methods introduced earlier on a set of realistic instances provided by EDF. \Cref{subsec:data} describes the available data and the procedure used to generate standardized instances from the raw EDF data, enabling a consistent evaluation of the methods. \Cref{subsec:acceleration} analyzes the impact of relaxing the dominance conditions and identifies the most effective variant of the CG-based algorithm. \Cref{subsec:performance} compares the resulting CG-based algorithm with the compact formulation solved by a commercial solver, considering both the weighted and sequential reformulations. Finally, \Cref{subsec:raw_instances} evaluates the performance of these four approaches on some of the largest real-life instances provided by EDF.

Since electricity utility companies typically generate technician schedules on a daily basis, computational efficiency is a key requirement. Accordingly, we consider a short 5-minute time limit, reflecting real-time operational constraints. For completeness, we also report results under an extended overnight time budget of 8 hours, although the proposed CG-based algorithm is designed for short runtimes and may not substantially benefit from longer runtimes. All computational experiments were carried out on a laptop equipped with an 11th Gen Intel\textsuperscript{\textregistered} Core\texttrademark{} i5-11500H processor (2.70~GHz, 6 cores) and 32~GB of RAM. The algorithms were implemented in Julia, and use Gurobi~v12 as the commercial MILP solver.

\subsection{Data} \label{subsec:data}

\begin{figure}[ht]
	\centering
	\begin{tikzpicture}
		\begin{axis}[
			width=14cm, height=8cm,
			xmin=0, xmax=300,
			ymin=0, ymax=60,
			xtick={0,50,100,150,200,250,300},
			ytick={0,10,20,30,40,50,60},
			grid=both,
			major grid style={gray!30},
			xlabel=\tiny{NUMBER OF INTERVENTIONS},
			ylabel=\tiny{NUMBER OF VEHICLES},
			legend style={
				at={(0.5,-0.18)},
				anchor=north,
				legend columns=3,
				draw=none,
				/tikz/column sep=8pt,   
				row sep=3pt,            
				inner xsep=4pt,         
				inner ysep=2pt,
			},
			legend cell align=left,
			]
			
			\pgfplotstableread[col sep=space]{
				I   V
				203 41
				184 43
				145 41
				141 27
				183 35
				139 37
				201 35
				243 47
				196 33
				205 37
				209 34
				182 41
				233 42
				161 33
				248 33
				252 43
				111 32
				150 31
				180 38
				224 43
				200 48
				231 45
				154 29
				77  20
				174 31
				194 37
				218 45
				171 34
				204 41
				202 40
				204 43
				204 39
				189 36
				161 36
				170 29
				215 39
				243 44
				169 39
				205 35
				97  34
			}\rawdata
			
			\addplot[
			only marks,
			mark=diamond*,
			mark size=2pt,
			blue
			] table[x=I,y=V]{\rawdata};
			\addlegendentry{\footnotesize{Raw Data}}
			
			\pgfplotstableread[col sep=space]{
				I   V
				20  4
				100 20
				150 30
				200 40
			}\testdata
			
			\addplot[
			only marks,
			mark=square*,
			mark size=2.5pt,
			red
			] table[x=I,y=V]{\testdata};
			\addlegendentry{\footnotesize{Test Instances}}
			
			\addplot[
			domain=0:300,
			samples=2,
			thick,
			orange
			] {x/5};
			\addlegendentry{\footnotesize{Linear (Ratio 1:5)}}
			
		\end{axis}
	\end{tikzpicture}
	\caption{Vehicle-to-intervention ratio of the raw and test (standardized) instances}
	\label{fig:raw_edf_data}
\end{figure}
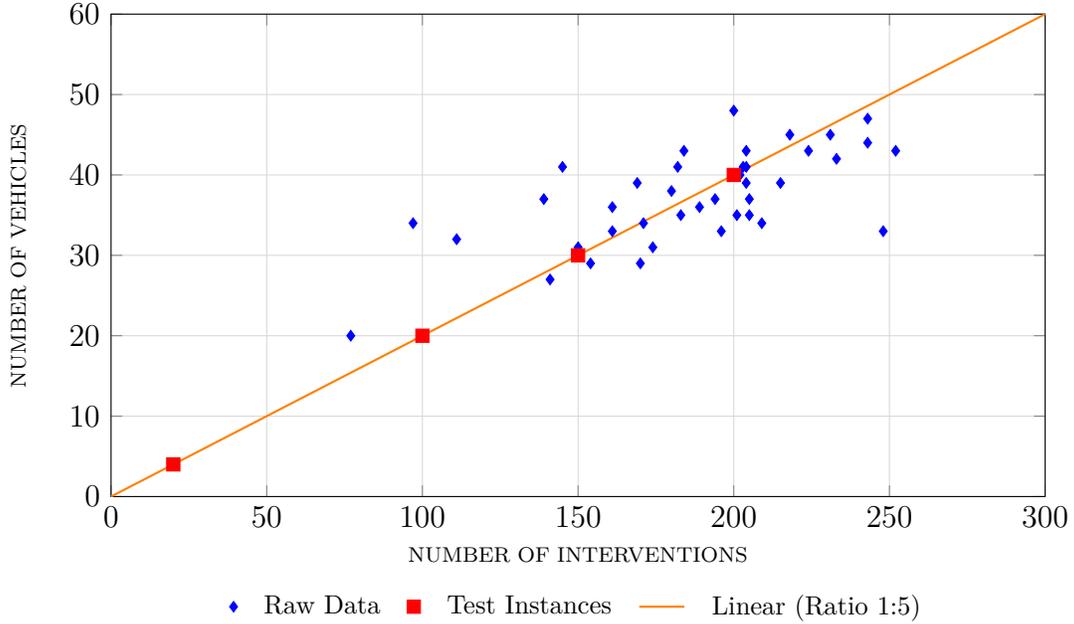

The number of daily interventions performed by EDF ranges from 77 to 252, while the number of vehicles varies between 20 and 48. From these data, we constructed a set of standardized test instances containing 20, 30, and 40 vehicles, with a fixed vehicle-to-intervention ratio of 1:5. To do so we first removed vehicles from the raw EDF instances to obtain the desired number. We then randomly removed interventions to precisely adjust the ratio to 1:5. For example, consider an original instance with 22 vehicles and 105 interventions. After removing two vehicles, the number of vehicles becomes 20, resulting in a ratio of 1:5.25. Subsequently, five interventions are removed to achieve the exact target ratio of 1:5. We then grouped the resulting instances into three categories according to their size: small (S; 100 interventions and 20 vehicles), medium (M; 150 interventions and 30 vehicles), and large (L; 200 interventions and 40 vehicles). We also generated extra-small instances (XS; 20 interventions and 4 vehicles), which are small enough to be solved to optimality and thus serve as benchmarks for evaluating solution quality.  Table~\ref{tab:instance_sizes} summarizes the main characteristics of these standardized instance categories. \Cref{fig:raw_edf_data} illustrates the vehicle-to-intervention ratio for the 40 raw instances provided by EDF and for the standardized test instances, highlighting the diversity of instance sizes (blue diamonds), the XS, S, M, and L instance categories (red rectangles) and the fixed 1:5 ratio (orange line). 

In addition to these controlled categories, the algorithm is tested on five of the largest EDF instances. These instances represent the upper bound on the problem size observed in practice and serve to assess the scalability of the proposed approach under extreme operational conditions. Each instance identifier encodes two key characteristics: the number of interventions and the number of vehicles involved, following the format: $Ix\_Vy,$ where $x$ denotes the number of interventions and $y$ the number of vehicles. For example, $I243\_V47$ corresponds to an instance with 243 interventions and 47 vehicles. The number of interventions and number of vehicles of the five evaluated instances range between 233 and 252, and 33 and 47, respectively.

\begin{table}[ht]
	\begin{center}
		\setlength{\tabcolsep}{4pt}
		\caption{Characteristics of the test instance categories.}
		\label{tab:instance_sizes}
		\begin{tabular*}{\textwidth}{@{\extracolsep{\fill}}lccccc}
			
			\toprule
			\textbf{Instance size} & \textbf{XS} &  \textbf{S} & \textbf{M}  & \textbf{L}  \\
			\midrule
			Number of interventions &  20 & 100 & 150 & 200 \\
			\midrule
			Number of vehicles & 4 & 20 & 30 & 40 \\ 
			\midrule
			Number test instances & 20 & 20 & 15 & 5 \\ 
			\bottomrule
		\end{tabular*}
	\end{center}
\end{table}

\subsection{Impact of the relaxation of dominance} \label{subsec:acceleration}

We first analyze the impact of relaxing dominance conditions (\cref{alg:dominance}) in the PS, which directly affects the size and diversity of the column pool explored by the CG-based algorithm. The comparison is conducted under a 5-minute time limit, evaluating the CG-based algorithm with both the weighted (CG W.) and sequential (CG S.) reformulations.

For each instance category, \cref{Table:relaxationch} reports the number of tested instances and the number of best-known solutions identified. The latter is expressed as a tuple $(a,b)$, where $a$ is the number of best solutions obtained for the primary objective ($f^1$) and $b$ is the number of best solutions obtained for the secondary objective ($f^2$) only when the best solution for the primary objective is obtained. We also report the mean gap relative to the best solution found within each instance category across all these four approaches. For the secondary objective, the mean gap is computed only when the best solution for the primary objective is available; otherwise, no value is reported (marked by --). Finally, the table also includes the average number of routes generated by the CG-based algorithm and the average number of iterations. 

\begin{table}[ht]
	\begin{center}
		\setlength{\tabcolsep}{4pt}
		\caption{Impact of the relaxation of dominance conditions on the CG Weighted and CG Sequential approaches.}
		\label{Table:relaxationch}
		\begin{tabular}{@{}clcccc@{}}
			\toprule
			\textbf{Inst.} (\#)&   \multicolumn{1}{c}{\textbf{Metric}} & 
			\multicolumn{2}{c}{\textbf{CG W.}} & \multicolumn{2}{c}{\textbf{CG S.}} \\
			\cmidrule(lr){3-4} \cmidrule(lr){5-6}
			& & \textbf{W/o Relax.} & \textbf{W/ Relax.} & \textbf{W/o Relax.} & \textbf{W/ Relax.} \\
			\midrule
			
			& Nb.\ of best solution & (20, 11) & (20, 15) & (18, 11) & (20, 17) \\
			\textbf{XS} (20) & Avg.\ best-gap (\%)   & (0.00, 2.34) & (0.00, 0.25) & (2.14, 0.53) & (0.00, 1.66) \\
			
			& Avg.\ routes & 1311 & 1480 & 1372 & 936 \\
			& Avg.\ iters. & 15 & 35 & 21 & 40 \\
			\midrule
			
			& Nb.\ of best solution & (0, 0) & (20, 9) & (1, 0) & (20, 11) \\
			\textbf{S} (20)  & Avg.\ best-gap (\%)       & (24.54, --) & (0.00, 1.74) & (24.25, 31.38) & (0.00, 2.18) \\
			& Avg.\ routes & 582 & 3690 & 558 & 3097 \\
			& Avg.\ iters. & 3 & 58 & 3 & 46 \\
			\midrule
			
			& Nb.\ of best solution & (0, 0) & (9, 5) & (0, 0) & (13, 10) \\
			{\textbf{M}}  (15) & Avg.\ best-gap (\%)       & (43.00, --) & (0.39, 1.21) & (41.85, --) & (0.25, 1.13) \\
			& Avg.\ routes & 359 & 6649 & 244 & 4269 \\
			& Avg.\ iters. & 1 & 57 & 1 & 37 \\
			\midrule
			
			& Nb.\ of best solution & (0, 0) & (3, 3) & (0, 0) & (2, 2) \\
			{\textbf{L}}   (5) & Avg.\ best-gap (\%)       & (16.37, --) & (0.22, 0.00) & (15.71, --) & (0.07, 0.00) \\
			& Avg.\ routes & 386 & 10650 & 246 & 5860 \\
			& Avg.\ iters. & 1 & 53 & 1 & 35 \\
			\bottomrule
		\end{tabular}
	\end{center}
\end{table}

The results show that applying \cref{alg:dominance} generally improves performance across most instance categories, although the effect remains limited for the XS instances. In this case, the search space is small and this fact restricts the number of feasible paths in general, reducing the potential impact of dominance relaxation. Indeed, the mean best-gap of CG W. decreases from $2.34\%$ to $0.25\%$, and the number of best solutions for the second objective increases from $11$ to $15$, confirming a modest but consistent improvement. In contrast, the effect is mixed in the CG S. variant: the first-objective gap decreases (from $2.14\%$ to $0.00\%$), but the second-objective gap slightly increases (from $0.53\%$ to $1.66\%$). Nevertheless, the number of best solutions for both objectives still improves, from $(18,11)$ to $(20,17)$.

For the S instances, the benefits become much clearer. CG W. reduces the mean best-gap for the first objective from $24.54\%$ to $0.00\%$, and the number of best solutions increases from $(0,0)$ to $(20,9)$. A similar pattern is observed for the CG S.: the mean best-gaps improve from $(24.25\%, 31.38\%)$ to $(0.00\%, 2.18\%)$, along with an increase in the number of best solutions from $(1,0)$ to $(20,11)$. The inferior performance observed when \cref{alg:dominance} is not applied can be attributed to the limited number of iterations and the restricted number of generated routes.

The impact of the algorithm on the M and L instances becomes even more apparent. When all dominance conditions are enforced from the beginning, the CG algorithm rarely discards paths, leading to long lists of unexpanded paths and almost no progress. As a result, CG W. performs only one iteration on average, generating 359 and 386 routes for the M and L instances, respectively, and CG S. generates similarly few routes (244 and 246). With the relaxation mechanism activated, the number of iterations increases to $57$ and $53$ for CG W., and to $37$ and $35$ for CG S. (for M and L instances, respectively). The number of generated routes also increases by an order of magnitude, giving the MP($\mathcal{R'}$) substantially richer column sets. This intensified route generation directly translates into superior solution quality. The mean best-gap decreases substantially: for CG W., from $43.00\%$ to $0.39\%$ on M instances and from $16.37\%$ to $0.22\%$ on L instances; for CG S., from $41.85\%$ to $0.25\%$ on M instances and from $15.71\%$ to $0.07\%$ on L instances.

Overall, relaxing the dominance conditions consistently improves the performance of both CG-based algorithm variants across all instance categories, allowing a more effective exploration of the solution space within the given time limit. Accordingly, this mechanism is retained in all subsequent experiments, as it provides the best trade-off between computational effort and solution quality.


\subsection{Performance} \label{subsec:performance}

\Cref{Table:comp_result_5min,Table:comp_result_8h} summarize the results of the comparison between the two variants of the compact formulation (Compact W. and Compact S., for the weighted and sequential reformulations, respectively) and the extended formulation (CG W. and CG S.) across the standardized instance categories, for both considered time limits.

\begin{table}[ht]
	\begin{center}
		\setlength{\tabcolsep}{4pt}
		\caption{Comparison of compact formulations and CG-based algorithm (extended formulations) for the XS, S, M, and L instances (5-minute time limit).}
		\label{Table:comp_result_5min}
		\begin{tabular}{@{}clcccc@{}}
			\toprule
			{\textbf{Inst.}} (\#) & \textbf{Metric} & \textbf{Compact W.} & \textbf{Compact S.} & \textbf{CG W.} & \textbf{CG S.} \\
			\midrule
			{\textbf{XS}} (20)  & Nb.\ of best solution & (20, 18) & (20, 15) & (20, 12) & (20, 14) \\
			& Avg.\ best-gap (\%)   & (0.00, 0.07) & (0.00, 1.37) & (0.00, 0.55) & (0.00, 1.96) \\
			\midrule
			{\textbf{S}}  (20)  & Nb.\ of best solution & (7, 0) & (1, 0) & (20, 9) & (20, 11) \\
			& Avg.\ best-gap (\%)   & (1.86, 24.14) & (2.82, 31.08) & (0.00, 1.74) & (0.00, 2.18) \\
			\midrule
			{\textbf{M}} (15)  & Nb.\ of best solution & (0, 0) & (0, 0) & (9, 5) & (13, 10) \\
			& Avg.\ best-gap (\%)   & (14.97, --) & (16.94, --) & (0.29, 1.21) & (0.12, 1.13) \\
			\midrule
			{\textbf{L}}  (5)   & Nb.\ of best solution & (0, 0) & (0, 0) & (3, 3) & (2, 2) \\
			& Avg.\ best-gap (\%)   & (38.62, --) & (78.13, --) & (0.76, 0.00) & (0.26, 0.00) \\
			\bottomrule
		\end{tabular}
	\end{center}
\end{table}

For XS instances, all variants achieve the best solution on all 20 instances for the first objective, showing essentially equivalent performance. Differences emerge on the second objective: Compact W. performs slightly better, attaining the best value on 18 out of 20 instances, compared with 15, 12, and 14 instances for Compact S., CG W., and CG S., respectively. Its mean best-gap is 0.07\%, versus 1.37\%, 0.55\%, and 1.96\% for the other three methods. It is worth mentioning that on those instances, optimality was proved only in 5 cases, both for the Compact W. and Compact S. methods. This number increases to 8 for CG W. and to 10 for CG S. Note that none of the larger instances was solved to optimality. 

For S instances, the behavior changes significantly. The compact formulations exhibit clear limitations: Compact W. attains the best first-objective value on only 7 out of 20 instances and never achieves the best second-objective value, while Compact S. performs even worse, with just one instance reaching the best first-objective value and none for the second. Their mean best-gaps on the second objective are substantial, at 24.14\% and 31.08\%, respectively. In contrast, both CG variants consistently attain the best first objective on all 20 instances and reach the best second objective on 9 (CG W.) and 11 (CG S.) instances, with mean best-gaps reduced to 1.74\% and 2.18\%. From S instances onward, CG clearly outperforms the compact formulations.

The gap between the approaches widens considerably for M and L instances. For both sizes, the compact formulations fail to reach any best-known solution on either objective and exhibit large average best-gaps: $14.97\%$ and $16.94\%$ for Compact W. and Compact S. on M instances, and even more pronounced best-gaps of $38.62\%$ and $78.13\%$ on L instances. Due to the lexicographic structure, no meaningful values are reported for the second objective. In contrast, the CG-based algorithm variants remain effective across both categories. For M instances, CG W. attains the best first-objective value on 9 out of 15 instances and the best second-objective value on 5 of them, with mean best-gaps of $0.29\%$ and $1.21\%$, while CG S. reaches 13 best first objective values and 10 best second-objective values, with best-gaps reduced to $0.12\%$ and $1.13\%$. The advantage of CG becomes even more pronounced for L instances: CG W. finds the best solution for both objectives on 3 out of 5 instances with mean best-gaps of $0.76\%$ and $0.00\%$, and CG S. achieves similarly strong performance, with best-gaps of only $0.26\%$ and $0.00\%$.

\begin{table}[ht]
	\begin{center}
		\setlength{\tabcolsep}{4pt}
		\caption{Comparison of compact formulations and CG-based algorithm (extended formulations) for the XS, S, M, and L instances (8-hour time limit).}
		\label{Table:comp_result_8h}
		\begin{tabular}{@{}clcccc@{}}
			\toprule
			{\textbf{Inst.}} (\#) & \textbf{Metric} & \textbf{Compact W.} & \textbf{Compact S.} & \textbf{CG W.} & \textbf{CG S.} \\
			\midrule
			{\textbf{XS}} (20) & Nb.\ of best solution & (20, 18) & (20, 16) & (20, 13) & (20, 11) \\
			& Avg.\ best-gap (\%)   & (0.00, <0.01) & (0.00, 1.33) & (0.00, 0.46) & (0.00, 0.75) \\
			\midrule
			{\textbf{S}}  (20) & Nb.\ of best solution & (17, 0) & (15, 0) & (20, 12) & (20, 8) \\
			& Avg.\ best-gap (\%)   & (0.16, 21.49) & (0.14, 20.73) & (0.00, 1.22) & (0.00, 0.93) \\
			\midrule
			{\textbf{M}} (15) & Nb.\ of best solution & (11, 0) & (6, 0) & (14, 4) & (14, 11) \\
			& Avg.\ best-gap (\%)   & (0.48, 17.80) & (0.66, 10.82) & (0.03, 7.07) & (0.11, 0.21) \\
			\midrule
			{\textbf{L}}  (5)  & Nb.\ of best solution  & (1, 0) & (0, 0) & (3, 2) & (3, 3) \\
			& Avg.\ best-gap (\%)   & (0.82, 9.54) & (2.60, --) & (0.13, 6.39) & (0.07, 0.00) \\
			\bottomrule
		\end{tabular}
	\end{center}
\end{table}

With an extended 8-hour time limit, the compact formulations show a substantial improvement on the primary objective. For XS instances, all four methods now systematically achieve the best first-objective values. The approaches exhibit only small residual best-gaps on the second objective—ranging from less than 0.01\% for Compact W. to 1.33\% for Compact S.--indicating that, given sufficient time, all methods are essentially equivalent on this smallest instance category. With this additional solution time, both compact approaches could solve 8 instances to optimality, whereas this number increases to 12 for CG W but remains at 10 for CG S.

On S instances, the compact formulations substantially improve on the first objective: Compact W.\ and Compact S.\ now find the best value on 17 and 15 instances, respectively, with mean best-gaps below $0.2\%$. However, they still do not attain any best solution for the second objective, and the mean best-gaps on it remain high ($21.49\%$ and $20.73\%$). The CG variants continue to dominate in terms of lexicographic quality: both CG W. and CG S. reach the best first-objective value on all 20 instances, and they also achieve the best second-objective value on 12 and 8 instances, with mean best-gaps for the second objective around $1\%$. Thus, even with an extended time budget, CG remains preferable for jointly optimizing both objectives on S instances.

For M and L instances, all methods improve compared with the 5-minute setting, but the compact and CG formulations continue to exhibit markedly different behaviors. On M instances, the compact formulations attain most of the best-known values for the primary objective. Compact W. and Compact S. attain the best first-objective value on 11 and 6 instances, with mean best-gaps of $0.48\%$ and $0.66\%$, respectively. However, they still perform poorly on the second objective, with mean best-gaps of $17.80\%$ and $10.82\%$, respectively. In contrast, CG W. and CG S. reach the best value on 14 of the 15 instances for the first objective, and CG S. also attains the best value on 11 instances for the second objective, with mean best-gaps reduced to $0.03\%$. A similar pattern is observed for L instances. Compact W. reaches the best value on only 1 of the 5 instances for the first objective (with a mean best-gap of $0.82\%$) and still shows a large gap on the second objective ($9.54\%$), while Compact S. fails to attain the best solution on any instance and maintains a non-negligible mean best-gap on the first objective with undefined values on the second. By comparison, CG W. and CG S. each achieve the best value on three instances for the primary objective and the best value on two and three instances for the second objective, respectively. 


\subsection{Large real-life instances} \label{subsec:raw_instances}

\begin{figure}[ht]
	\centering
	\begin{tikzpicture}
		\begin{axis}[
			ybar,
			bar width=3pt,
			width=\linewidth,
			height=7cm,
			xlabel={\scriptsize Large EDF instances},
			ylabel={\small Objective value $f^1$},
			symbolic x coords={I243\_V47,I233\_V42,I248\_V33,I252\_V43,I243\_V44},
			xtick=data,
			ymin=0, ymax=10000,
			scaled y ticks=false,
			enlarge x limits=0.2,
			tick label style={font=\tiny},
			legend style={font=\tiny, at={(0.5,-0.25)}, anchor=north, legend columns=2},
			grid=major,
			bar shift auto
			]
			
			\addplot+[ybar, bar width=3pt, draw=CompactWColor, fill=CompactWColor] coordinates {
				(I243\_V47,2828)
				(I233\_V42,5653)
				(I248\_V33,4928)
				(I252\_V43,3206)
				(I243\_V44,5072)
			};
			\addlegendentry{Compact W.}
			
			\addplot+[ybar, bar width=3pt, draw=CompactSColor, fill=CompactSColor] coordinates {
				(I243\_V47,0)           
				(I233\_V42,3973)
				(I248\_V33,180)
				(I252\_V43,0)           
				(I243\_V44,60)
			};
			\addlegendentry{Compact S.}
			
			\addplot+[ybar, bar width=3pt, draw=CGWColor, fill=CGWColor] coordinates {
				(I243\_V47,8929)
				(I233\_V42,7801)
				(I248\_V33,7480)
				(I252\_V43,9425)
				(I243\_V44,9206)
			};
			\addlegendentry{CG W.}
			
			\addplot+[ybar, bar width=3pt, draw=CGSColor, fill=CGSColor] coordinates {
				(I243\_V47,9084)
				(I233\_V42,7986)
				(I248\_V33,8172)
				(I252\_V43,9540)
				(I243\_V44,9306)
			};
			\addlegendentry{CG S.}
			
		\end{axis}
	\end{tikzpicture}
	\caption{Objective values $f^1$ obtained on the five largest EDF instances (5-minute time limit).}
	\label{fig:raw_5min}
\end{figure}
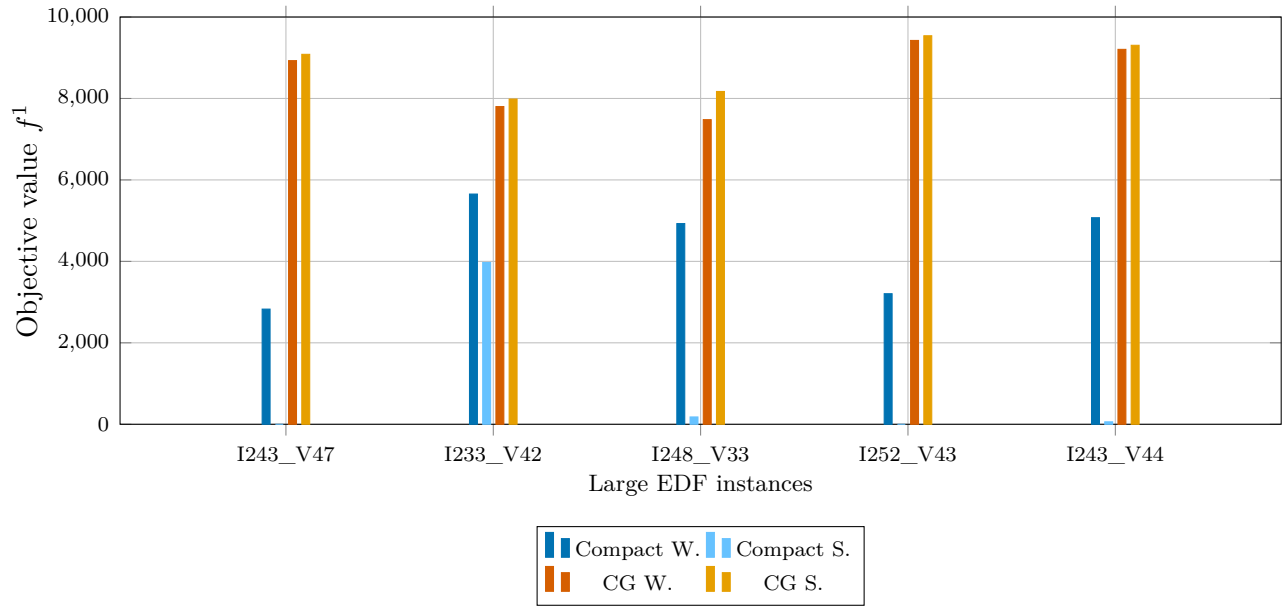

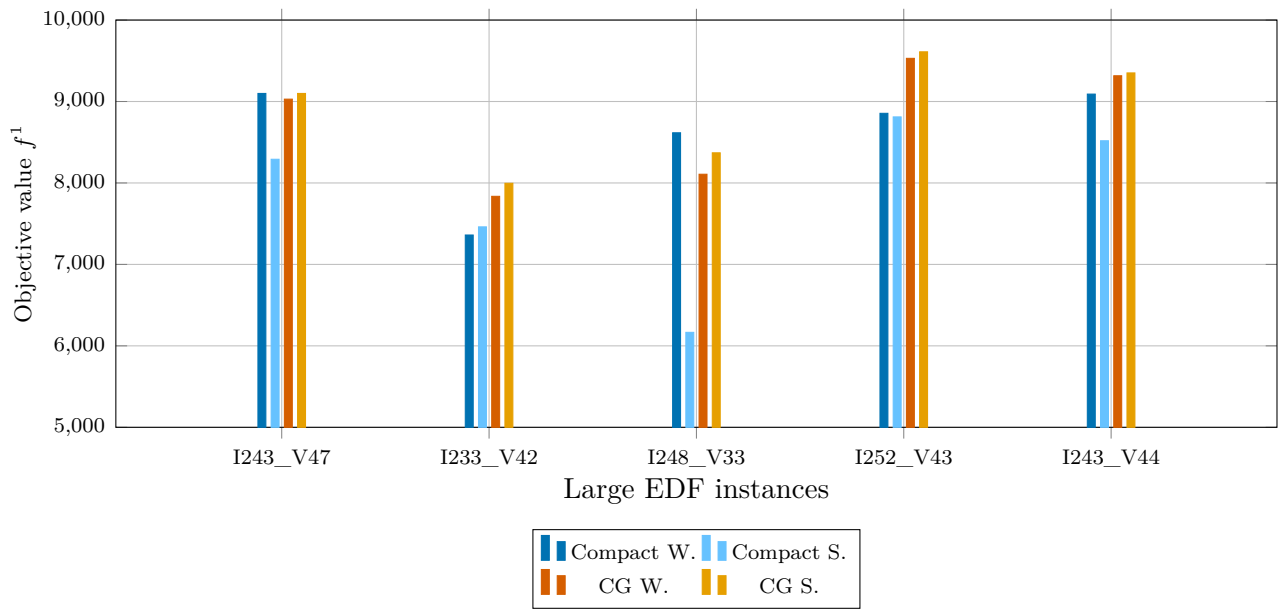
\begin{figure}[ht]
	\centering
	\begin{tikzpicture}
		\begin{axis}[
			ybar,
			bar width=3pt,
			width=\linewidth,
			height=7cm,
			xlabel={\small Large EDF instances},
			ylabel={\scriptsize Objective value $f^1$},
			symbolic x coords={I243\_V47,I233\_V42,I248\_V33,I252\_V43,I243\_V44},
			xtick=data,
			ymin=5000, ymax=10000,
			scaled y ticks=false,
			enlarge x limits=0.2,
			tick label style={font=\tiny},
			legend style={font=\tiny, at={(0.5,-0.25)}, anchor=north, legend columns=2},
			grid=major,
			bar shift auto
			]
			
			\addplot+[ybar, bar width=3pt, draw=CompactWColor, fill=CompactWColor] coordinates {
				(I243\_V47,9099)
				(I233\_V42,7361)
				(I248\_V33,8616)
				(I252\_V43,8855)
				(I243\_V44,9091)
			};
			\addlegendentry{Compact W.}
			
			\addplot+[ybar, bar width=3pt, draw=CompactSColor, fill=CompactSColor] coordinates {
				(I243\_V47,8291)
				(I233\_V42,7461)
				(I248\_V33,6166)
				(I252\_V43,8811)
				(I243\_V44,8518)
			};
			\addlegendentry{Compact S.}
			
			\addplot+[ybar, bar width=3pt, draw=CGWColor, fill=CGWColor] coordinates {
				(I243\_V47,9029)
				(I233\_V42,7836)
				(I248\_V33,8107)
				(I252\_V43,9530)
				(I243\_V44,9316)
			};
			\addlegendentry{CG W.}
			
			\addplot+[ybar, bar width=3pt, draw=CGSColor, fill=CGSColor] coordinates {
				(I243\_V47,9099)  
				(I233\_V42,7996)
				(I248\_V33,8370)
				(I252\_V43,9610)
				(I243\_V44,9351)
			};
			\addlegendentry{CG S.}
			
		\end{axis}
	\end{tikzpicture}
	\caption{Objective values $f^1$ obtained on the five largest EDF instances (8-hour time limit). The $y$-axis has been truncated to improve visibility of differences.}
	\label{fig:raw_8h}
\end{figure}

We next evaluate the performance of the four methods on the five largest EDF instances introduced in \Cref{subsec:data}. These experiments place particular emphasis on the 5-minute time limit, which reflects the tight computational requirements of operational planning. For completeness, we also report results obtained under an extended overnight time budget of 8 hours. Figures~\ref{fig:raw_5min} and~\ref{fig:raw_8h} present the corresponding values of the primary objective~$f^1$.

Under the 5-minute time limit (\Cref{fig:raw_5min}), the compact formulations show clear limitations on these instances. Compact~W.\ consistently produces lower values of $f^1$ than the CG-based methods, while Compact~S.\ sometimes returns solutions of very limited quality. This behavior for Compact~S.\ can be partly attributed to structural properties of the formulation. In particular, the formulation exhibits a high degree of symmetry, leading to a large number of equivalent solutions that are difficult to distinguish during the branch-and-bound process. This can slow down convergence and hinder the rapid identification of good feasible solutions. As a result, Compact~S.\ appears less competitive under very short computation times on large instances.

In contrast, both CG variants rapidly identify high-quality solutions: for all five instances, CG W.\ and CG S.\ deliver objective values that are clearly superior to those of the compact formulations, and their respective values remain close across instances. This behavior can be compared to the results observed for the L instances in \Cref{Table:comp_result_5min} and confirms that, under tight real-time constraints, the CG-based algorithm scales much more effectively than the compact formulation.

When the time limit is extended to 8 hours (\Cref{fig:raw_8h}; note that the scale of the objective value $f^1$ starts at 5000 for a clearer visualization of the differences for higher objective values), the relative advantage of the CG-based approaches becomes less pronounced compared with the 5-minute setting. The compact formulations, particularly the weighted variant, benefit more substantially from the additional runtime and are able to close a significant portion of the gap on the primary objective, with deviations remaining below 9\% in all cases. Compact S. also improves but remains less competitive, with deviations typically ranging from $6\%$ to $28\%$ across the five instances. As a result, the performance differences between compact and CG-based methods are smaller than under strict operational constraints.

This behavior is consistent with the design philosophy of the proposed CG-based algorithm, which was primarily developed and tuned for short computational horizons. Under such conditions, it rapidly generates high-quality columns and delivers strong feasible solutions. When longer runtimes are available, the compact formulations are given sufficient time to overcome their initial limitations, thereby reducing the observable performance gap.

Nevertheless, both CG variants continue to provide good results on the considered instances: CG~W.\ attains objective values within approximately $2\%$ of the best solution on every instance, while CG~S.\ reaches the best $f^1$ value on four of the five instances.  Increasing the time limit from 5 minutes to 8 hours yields only marginal additional improvements for both CG~W.\ and CG~S.\, suggesting that most algorithmic gains are already achieved within short computation times. 

In some cases, different methods attain identical values for the primary objective, making the secondary objective relevant for comparison. This occurs for instance $I243\_V47$ under the 8-hour time limit, where both Compact~W.\ and CG~S.\ reach the same best first-objective value $f^1 = 9099$ (\Cref{fig:raw_8h}). While Compact~W.\ yields $f^2 = 30\,994.5$, CG~S.\ achieves a lower cost value of $f^2 = 30\,261.5$, corresponding to an improvement of approximately $2.4\%$. 

In summary, the experiments on this set of real-life instances corroborate and extend the conclusions obtained from the standardized categories. Under a strict 5-minute time limit, the CG-based algorithm variants consistently produce high-quality solutions on the largest instances. When the computation time is increased to 8 hours, the relative advantage of CG progressively diminishes, as the compact formulations are able to further improve their solutions and reach objective values closer to the best-known ones.

\subsection{Managerial Insights}

This section presents managerial insights derived from the solutions obtained with the column generation approach based on the sequential reformulation (CG~S.). This methodological choice is supported by the computational results reported in \cref{subsec:raw_instances}, which indicate that CG~S.\ consistently achieves the highest number of best solutions compared to the other approaches under the same computational time limits. The analysis is conducted on the raw EDF instances, and the objective is to interpret, in operational and managerial terms, what is gained by allowing longer computational times on large real-life instances. 

For each raw EDF instance, we compare two solutions: $S_{5m}$, obtained under the short operational time limit of five minutes, and $S_{8h}$, obtained under the extended overnight run of eight hours. All numerical results referenced below are reported in Table~\ref{tab:MI_per_instance}. The problem having a lexicographic objective function, the primary objective is to maximize the executed intervention time, while the secondary objective is to minimize total operational cost. Hence, from a managerial perspective, improvements in solution quality therefore correspond to (i) an increase in executed intervention time or, given equal primary objective values, (ii) a reduction in total operational cost. The reported difference is defined as $\Delta = S_{8h} - S_{5m}$.

\begin{table}[ht]
	\begin{center}
		\setlength{\tabcolsep}{4pt}
		\caption{Comparison of CG~S.\ solutions after 5 minutes ($S_{5m}$) and 8 hours ($S_{8h}$) on raw EDF instances.}
		\label{tab:MI_per_instance}
		\begin{tabular}{@{}clccc@{}}
			\toprule
			\textbf{Instance} & \textbf{Metric} & $S_{5m}$ & $S_{8h}$ & $\Delta = S_{8h} - S_{5m}$ \\
			\midrule
			
			\multirow{2}{*}{\textbf{I243\_V47}}
			& Executed intervention time (min) & 9084 & 9099 & $+15$ \\
			& Total operational cost (\euro)   & 33759 & 30261.5 & $-3497.5$ \\
			\midrule
			
			\multirow{2}{*}{\textbf{I233\_V42}}
			& Executed intervention time (min) & 7986 & 7996 & $+10$ \\
			& Total operational cost (\euro)   & 31704 & 25585 & $-6119$ \\
			\midrule
			
			\multirow{2}{*}{\textbf{I248\_V33}}
			& Executed intervention time (min) & 8172 & 8370 & $+198$ \\
			& Total operational cost (\euro)   & 27569 & 27655.5 & $+86,5$ \\
			\midrule
			
			\multirow{2}{*}{\textbf{I252\_V43}}
			& Executed intervention time (min) & 9540 & 9610 & $+70$ \\
			& Total operational cost (\euro)   & 34952 & 35794 & $+842$ \\
			\midrule
			
			\multirow{2}{*}{\textbf{I243\_V44}}
			& Executed intervention time (min) & 9306 & 9351 & $+45$ \\
			& Total operational cost (\euro)   & 35623.5 & 34409.5 & $-1214$ \\
			\bottomrule
		\end{tabular}
	\end{center}
\end{table}

The observed differences between $S_{5m}$ and $S_{8h}$ fall into three distinct solution patterns:
\begin{itemize}
	\item The first pattern is illustrated by instance I248\_V33, for which extended computation time enables the recovery of a substantial amount of additional intervention time at a marginal increase in cost. Compared with $S_{5m}$, $S_{8h}$ achieves an additional 2.4\% (198 minutes) of executed intervention time, while the total operational cost increases by 0.3\% (86,5\euro).
	This result suggests that, for certain instances, a short computational time limit may leave a non-negligible portion of feasible interventions unassigned.
	
	\item The second pattern concerns instances in which the increase in executed intervention time from $S_{5m}$ to $S_{8h}$ is comparatively small—typically between $0.1\%$ and $1\%$—yet the reduction in total operational cost is substantial. This behavior is observed for instances I243\_V47, I233\_V42, and I243\_V44, where cost reductions range from approximately $3\%$ to more than $20\%$. In these cases, additional computational time mainly enables the algorithm to reorganize routes and team assignments, yielding significantly lower costs while maintaining nearly the same level of service coverage.
	
	\item The third pattern is exemplified by instance I252\_V43, in which the algorithm achieves a marginal increase in executed intervention time (+0.73\%) at the expense of a more substantial rise in total operational cost (+2.40\%). Although the improvement in service time is limited, this outcome remains fully consistent with the lexicographic structure of the problem: any enhancement of the primary objective—maximizing executed intervention time—takes precedence over the secondary objective of cost minimization. Consequently, an increase in cost is accepted whenever it enables additional interventions to be completed. From a managerial perspective, this result highlights the importance of coherence between the stated strategic priorities and the implemented decisions. If maximizing service coverage is indeed the primary objective, then accepting a higher operational cost—even for a modest gain in executed intervention time—is consistent with the chosen decision framework.
\end{itemize}

Collectively, these three solution patterns highlight the algorithm’s ability to reconfigure team allocations and routing decisions when additional computational time is available. They also indicate that, although the CG-based approach already performs strongly within five minutes, additional computational time can still yield non-negligible improvements on certain large instances, either in terms of increased intervention coverage or enhanced cost efficiency.

\section{Conclusion and future work} \label{sec:conclusion}
In this paper, we study a multi-objective variant of the TRSP in which the objectives are optimized in lexicographic order. The problem incorporates operational constraints reflecting the daily technical interventions performed by EDF, the French electric utility company. We introduce both a compact and an extended MILP formulation, the latter relying on the set-packing formulation commonly used in routing problems. For each formulation, we derive weighted and sequential single-objective reformulations that preserve the lexicographic priority of the objectives. The extended reformulations are solved using an exact CG–based algorithm. The proposed methods are evaluated on a set of standardized instances generated from real EDF data, and their scalability and robustness are further assessed on the five largest real-life instances provided by EDF.

A primary objective of this work was to design an approach capable of delivering high-quality solutions under tight operational constraints. Under the 5-minute time limit, the computational experiments presented in \Cref{subsec:performance} show that, for the smallest instances (XS), all solution methods perform similarly on the primary objective, with differences emerging only on the secondary objective. Among the compact formulations, the weighted reformulation consistently outperforms the sequential one, a trend that also holds for larger instance categories. Indeed, the weighted reformulation attains the best first-objective value on more instances and exhibits lower average best-gaps. For the S, M, and L instance categories, the CG–based algorithm variants clearly outperform their compact counterparts. More precisely, the sequential reformulation proves more effective overall, as it consistently finds the best-known solutions on more instances and achieves lower mean best-gaps.

When the time budget is extended to eight hours, the performance gap between the compact and CG-based approaches becomes less pronounced. As in the 5-minute setting, the compact weighted reformulation outperforms the sequential one across all instance categories. Nevertheless, both CG-based variants exhibit very similar performance levels, with the weighted reformulation yielding slightly better results for the XS and S instances and the sequential reformulation for the M and L ones. On the large real-life instances, the compact formulations are able to leverage the longer runtime to further improve their solutions, thereby reducing the observable advantage of CG. This outcome is consistent with the design of the proposed algorithm, which was primarily developed and tuned for short runtimes, rather than for sustained long-run convergence.

The managerial analysis of the large real-life instances reveals that extending computation from 5 minutes to 8 hours can still generate non-negligible improvements on certain instances, both in terms of executed intervention time and operational cost.
These findings suggest that, although the CG-based approach already captures most structural gains rapidly, additional optimization potential remains when longer computation times are available.

This observation motivates future research aimed at improving the long-run performance of the CG-based algorithm. In extended runs, the number of generated columns increases substantially, which may complicate the resolution of the integer master problem and slow down convergence. Integrating advanced mechanisms—such as diving heuristics, tailored branching strategies, or interleaved column generation and integer re-optimization phases—could further improve solution quality when larger time budgets are available.

An alternative to the single-objective reformulations considered in this paper is the development of a fully lexicographic CG approach inspired by \cite{tellache2024linear}. As described in Section~\ref{sec:LR}, the authors introduce an exact method for integer lexicographic linear programming problems. In this context, the continuous relaxation is formulated as a linear lexicographic programming problem, and the pricing subproblem is modeled as a resource-constrained longest path problem with lexicographic costs. The method is applied to the airline crew scheduling problem under a preferential bidding system and demonstrates competitive performance with state-of-the-art techniques on industrial instances ranging from 17 to 150 pilots (objectives). However, given that our problem involves only two objectives, the practical advantage of a fully lexicographic CG approach may be restricted compared with the single-objective reformulations considered here.

Finally, an additional aspect concerns the treatment of heterogeneity and fairness in workload distribution across vehicles and technicians. From a managerial perspective, disparities in workload allocation may be undesirable, as they can compromise perceived fairness, schedule acceptability, and long-term operational efficiency. An interesting extension of this work would therefore be to incorporate fairness-related criteria into the optimization process, enabling decision-makers to explicitly balance efficiency and equity in the resulting schedules.

\section*{Acknowledgments}

This research was supported by the Gaspard Monge Program for Optimization, Operations Research and their Interactions with Data Science (PGMO) of the Fondation Mathématique Jacques Hadamard (FMJH). The authors also gratefully acknowledge the OSIRIS Department at EDF for the valuable discussions conducted throughout the project and for providing the data used in this study.
	
	\bibliographystyle{unsrtnat}
	\bibliography{references}	
\end{document}